\documentclass[11pt]{article}

\usepackage[T1]{fontenc}
\usepackage{graphicx}
\usepackage{amsmath}
\usepackage{amsfonts}
\usepackage{amsthm}
\usepackage{amssymb}
\usepackage{amscd}
\usepackage{epsfig}
\usepackage{verbatim}
\usepackage{fancybox}
\usepackage{moreverb}
\usepackage{psfrag}
\usepackage{cellspace}
\usepackage{hyperref}
\usepackage{latexsym}
\usepackage{psfrag}
\usepackage{tikz}
\usetikzlibrary{arrows,positioning} 
\tikzset{
    >=stealth',
    punkt/.style={
           rectangle,
           rounded corners,
           draw=black, very thick,
           text width=4.9cm,
           minimum height=1cm,
           text centered},
    pil/.style={
           ->,
           thick,
           shorten <=2pt,
           shorten >=2pt,}
}
\tikzset{
    >=stealth',
    punktt/.style={
           rectangle,
           rounded corners,
           draw=black, very thick,
           text width=1cm,
           minimum height=1cm,
           text centered},
    pill/.style={
           ->,
           thick,
           shorten <=2pt,
           shorten >=2pt,}
}

\definecolor{blanc}{rgb}{1.,1.,1.}
\definecolor{vert}{rgb}{0.,0.5,0.}
\definecolor{rouge}{rgb}{0.8,0.,0.}
\definecolor{violet}{rgb}{0.5,0.,0.4}
\definecolor{bleu}{rgb}{0.,0.,0.5}
\definecolor{orange}{rgb}{0.8,0.4,0.}
\definecolor{light-blue}{rgb}{0.5,0.5,0.7}
\definecolor{light-red}{rgb}{0.8,0.4,0.4}
\definecolor{noir}{rgb}{0.,0.,0.}
\definecolor{gris}{rgb}{0.8.,0.8,0.7}

\newcommand{\eps}{\varepsilon}
 
\newcommand{\Avec}{\mathbf A} 
 
\newcommand{\Bvec}{\mathbf B}

\newcommand{\Evec}{\mathbf E} 
\newcommand{\fvec}{\mathbf f} 
\newcommand{\Fvec}{\mathbf F} 
 
\newcommand{\Ivec}{\mathbf I} 
\newcommand{\pvec}{\mathbf p} 
\newcommand{\Pvec}{\mathbf P} 
\newcommand{\ptrf}{\boldsymbol \pi} 
\newcommand{\qvec}{\mathbf q} 
\newcommand{\Qvec}{\mathbf Q} 
\newcommand{\rvec}{\mathbf r} 
\newcommand{\rtrf}{\boldsymbol\rho} 
\newcommand{\Rvec}{\mathbf R} 
\newcommand{\xvec}{\mathbf x} 
\newcommand{\Xvec}{\mathbf X} 
\newcommand{\xtrf}{\boldsymbol\xi} 
\newcommand{\xtrfb}{\boldsymbol\upsilon} 
\newcommand{\vvec}{\mathbf v} 
\newcommand{\Vvec}{\mathbf V} 
\newcommand{\yvec}{\mathbf y} 
 
\newcommand{\ytrf}{\text{\small$\boldsymbol\Upsilon$}} 
\newcommand{\ytrfb}{\text{\small$\boldsymbol\lambda$}} 
\newcommand{\zvec}{\mathbf z} 
\newcommand{\Zvec}{\mathbf Z} 
\newcommand{\ztrf}{\boldsymbol\zeta}

\newcommand{\zetavec}{\boldsymbol\zeta}

\newcommand{\Ocal}{{\cal O}}
\newcommand{\Pcal}{{\cal P}}

\newcommand{\Lcal}{{\cal L}}
\newcommand{\Tcal}{{\cal T}}
\newcommand{\Vcal}{{\cal V}}
\newcommand{\Scal}{{\cal S}}
\newcommand{\Xcal}{{\cal X}}
\newcommand{\Xcalvec}{{\boldsymbol{\cal X}}}

\newcommand{\rit}{\mathbb R} 
 
\newcommand{\Nabla}{\nabla\hspace{-2pt}}

\newcommand{\ds}{\displaystyle}
\newcommand{\fracp}[2]{\frac{\partial #1}{\partial #2}}
\newcommand{\poibrack}[2]{{\left\{{#1},{#2}\right\}}}

\newtheorem{remark}{Remark}[section]

\newtheorem{theorem}{Theorem}[section]

\begin{document}

\title{The Gyro-Kinetic Approximation\\{\Large An attempt at explaining the method based on Darboux Algorithm and Lie Transform}}
\author{Emmanuel Frénod\thanks{Universit\'e Europ\'eenne de Bretagne, Lab-STICC (UMR CNRS 3192), 
Universit\'e de Bretagne-Sud, Centre Yves Coppens, Campus de Tohannic, F-56017, Vannes \& 
Projet INRIA Calvi, Universit\'{e} de Strasbourg, IRMA, 7 rue Ren\'e Descartes, F-67084 Strasbourg Cedex, France.} 
\and 
Mathieu Lutz\thanks{Universit\'{e} de Strasbourg, IRMA,  7 rue Ren\'e Descartes, F-67084 Strasbourg Cedex, France \& Projet INRIA Calvi.}}
\date{}
\maketitle
\begin{center}{\it Summer school Fusion - Paris - September 2011} \end{center}

~

\abstract{\scriptsize  This Proceeding presents the method that allows us to get the Gyro-Kinetic Approximation of the Dynamical System
satisfied by the trajectory of a particle submitted to a Strong Magnetic Field. The goal of the method is to build a change of coordinates in order
to make the dynamic of two components of the trajectory to disappear. This change of coordinates is based on a Darboux mathematical Algorithm
and on a Lie Transform} 

%
%
\section{Introduction}

\paragraph{Scientific framework -} At the end of the 70',
{Littlejohn \cite{littlejohn:1979,littlejohn:1981,littlejohn:1982}}
shed new light on  what is called \emph{the Guiding Center Approximation.}
His approach was based on mathematical theories - Hamiltonian Mechanics, Differential Geometry, Symplectic Geometry - 
in order to clarify what has been done for years in the domain 
(see {Kruskal \cite{Kruskal1965}},
{Gardner \cite{Gardner1959}},
{Northrop  \cite{northrop:1961}},
{Northrop \& Rome \cite{10.1063/1.862226}}).
His papers claim that it has been for him an enormous effort to reach this goal, since he had to incorporate into a physical affordable
theory high level mathematics.
Sure, this theory is a nice success. It has been being  widely used by physicists to deduce related models 
(\emph{Finite Larmor Radius Approximation, Drift Kinetic Model, Quasi-Neutral Gyro-kinetic Model, etc.}, see for instance 
{Koseleff \cite{Koseleff94comparisonbetween}},
{Brizard \cite{10.1063/1.871465}},
{Dubin \emph{et al.} \cite{dubin/etal:1983}},
{Frieman  \& Chen \cite{10.1063/1.863762}},
{Hahm \cite{10.1063/1.866544}},
{Hahm, Lee \& Brizard \cite{10.1063/1.866641}},
{Parra \& Catto \cite{0741-3335-50-6-065014,0741-3335-51-6-065002,0741-3335-52-4-045004}})
 making up the \emph{Gyro-Kinetic Approximation Theory}, which is the basis of
all kinetic codes used to simulate Plasma Turbulence emergence and evolution in Tokamak
(see for instance {Grandgirard \emph{et al.} \cite{Grandgirard2006395,0741-3335-49-12B-S16}}).

The incorporation of mathematical concepts into a physical theory has been done so nicely that the
resulting Gyro-Kinetic Approximation Theory is now very difficult for mathematicians.

\begin{figure}[htbp]
\begin{center}
\includegraphics[width=4cm] {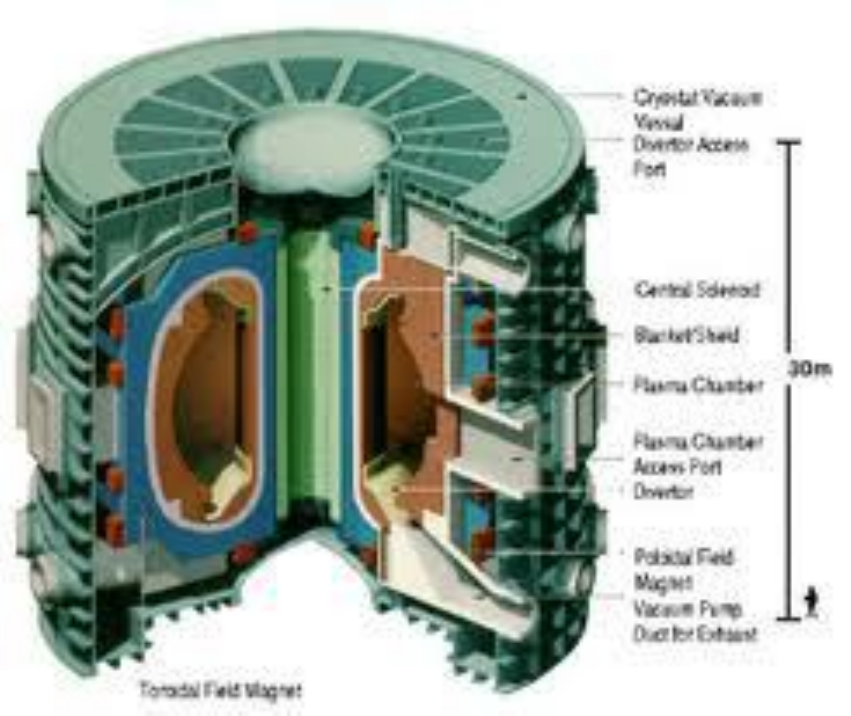}
\caption{Artist vision of Iter.}
\label{figter4}
\end{center}
\end{figure}

This talk (and proceeding) is the first step into the rewriting of the Gyro-Kinetic Approximation Theory
into a mathematically affordable theory.
It is only a summarize, to understand, with a mathematical slant, several aspects of what is done in the references just cited.
%
\paragraph{Charge particles submitted to Strong Magnetic Field -} \label{par:cpssmf}
The context of the Gyro-kinetic Approximation is Tokamak Physics. An artist vision of Iter, which is a Tokamak, is given in
Figure \ref{figter4}. The vessel of a Tokamak is the interior of a torus with a vertical axis of symmetry. Along the 
torus, electromagnets can generate a large magnetic field.
\label{dede}\\
We begin by considering a charged particle within the vessel of a Tokamak in service.
In Usual Coordinates $(\xvec,\vvec) = (x_1,x_2,x_3, v_1,v_2,v_3)$, where $\xvec$ stands for
the position variable and $\vvec$ for the velocity variable, the position and the velocity at time $t$
$(\ds \Xvec(t;\xvec,\vvec,s), $ $\Vvec(t;\xvec,\vvec,s))$ of the particle which is in $\xvec$ with velocity $\vvec$ at time $s$
is the solution to
 \begin{align} 
     &\fracp{\Xvec}{t}=\Vvec,
     \label{DsDim1} 
     \\
     &\fracp{\Vvec}{t}= \frac qm (\Evec(\Xvec)+\Vvec \times \Bvec(\Xvec)).
     \label{DsDim2} 
\end{align}
This dynamical system means nothing but that time derivative of the position is the velocity and that  time derivative of the velocity, 
which is the acceleration,
in linked with Lorentz force by Newton's law.
In this dynamical system magnetic field $\Bvec$  is composed of a strong applied piece, of a strong self induced piece 
and of a self induced perturbations. 
The strong applied piece is generated by the electromagnets and wind around the axis of symmetry of the tokamak.
Charged particles winding, with a relatively large velocity, around the axis of symmetry within the vessel generate
the indispensable vertical component of the large magnetic field. This is what is called here the strong self induced piece.
The part of the self induced magnetic field which is not vertical is called here the self induced perturbations and is, from now, 
forgotten.
Electric field $\Evec$ is induced by the particles of the vessel.

Without loss of generality, we can consider that $\Evec$ is the opposite of the gradient of an electric potential, and $\Bvec$
the curl of vector potential, i.e.:
\begin{gather}
\label{PotDim} 
\Evec = -\nabla \Phi,  ~~\Bvec = \nabla\times\Avec.
\end{gather}
\paragraph{Helicoidal trajectories - Larmor Radius -}
It is well  known that the trajectory of a charged particle submitted to a magnetic field is a helix.
The helix axis is
the magnetic field direction and its radius, which is called in the context of Tokamak Plasma the Larmor Radius,
equals the norm of the projection of the velocity on the plan orthogonal to the magnetic field divided by the particle's mass
time the norm of the magnetic field. 
\begin{figure}[htbp]
\begin{center}
\includegraphics[width=6cm]{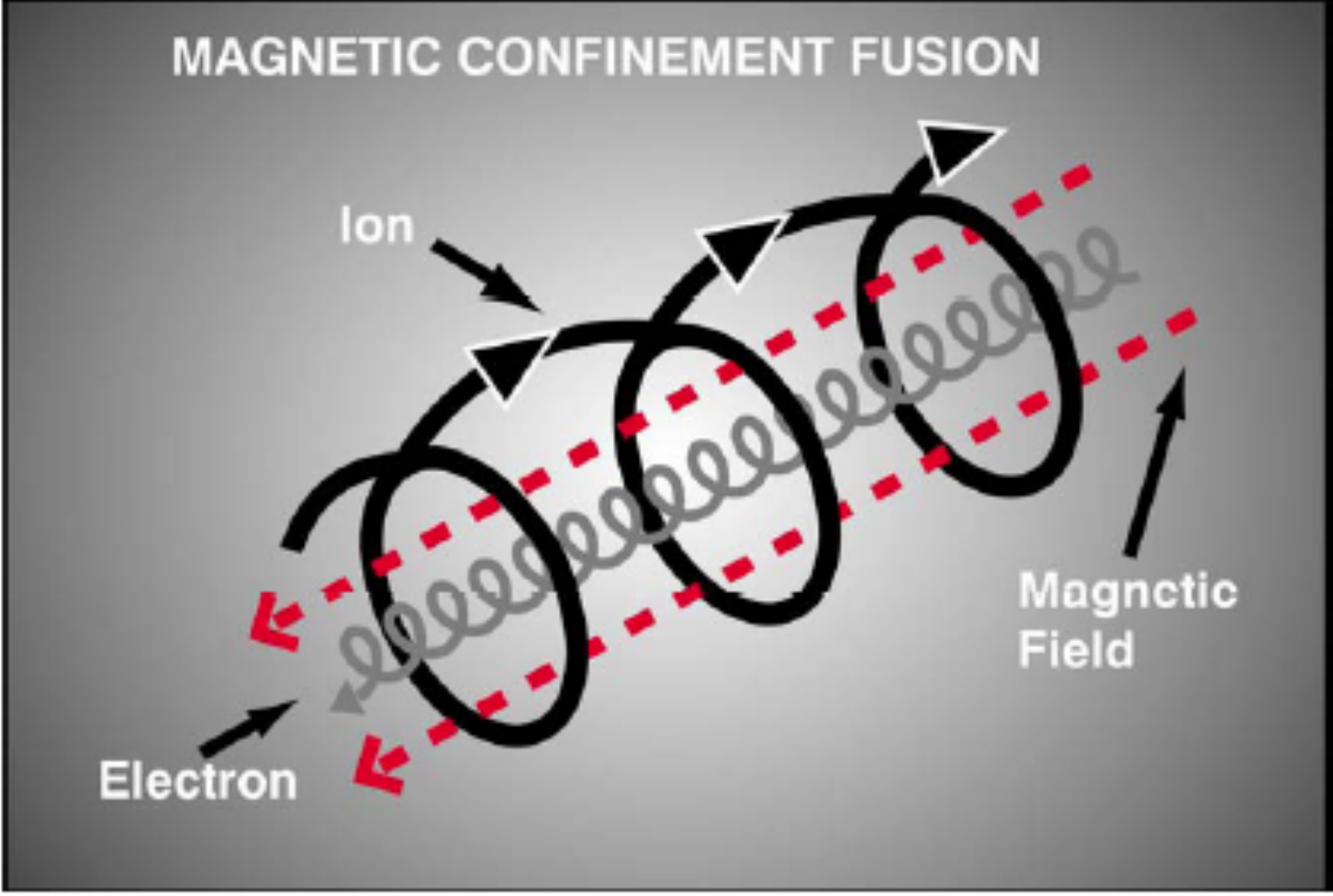}~~
\includegraphics[width=6cm]{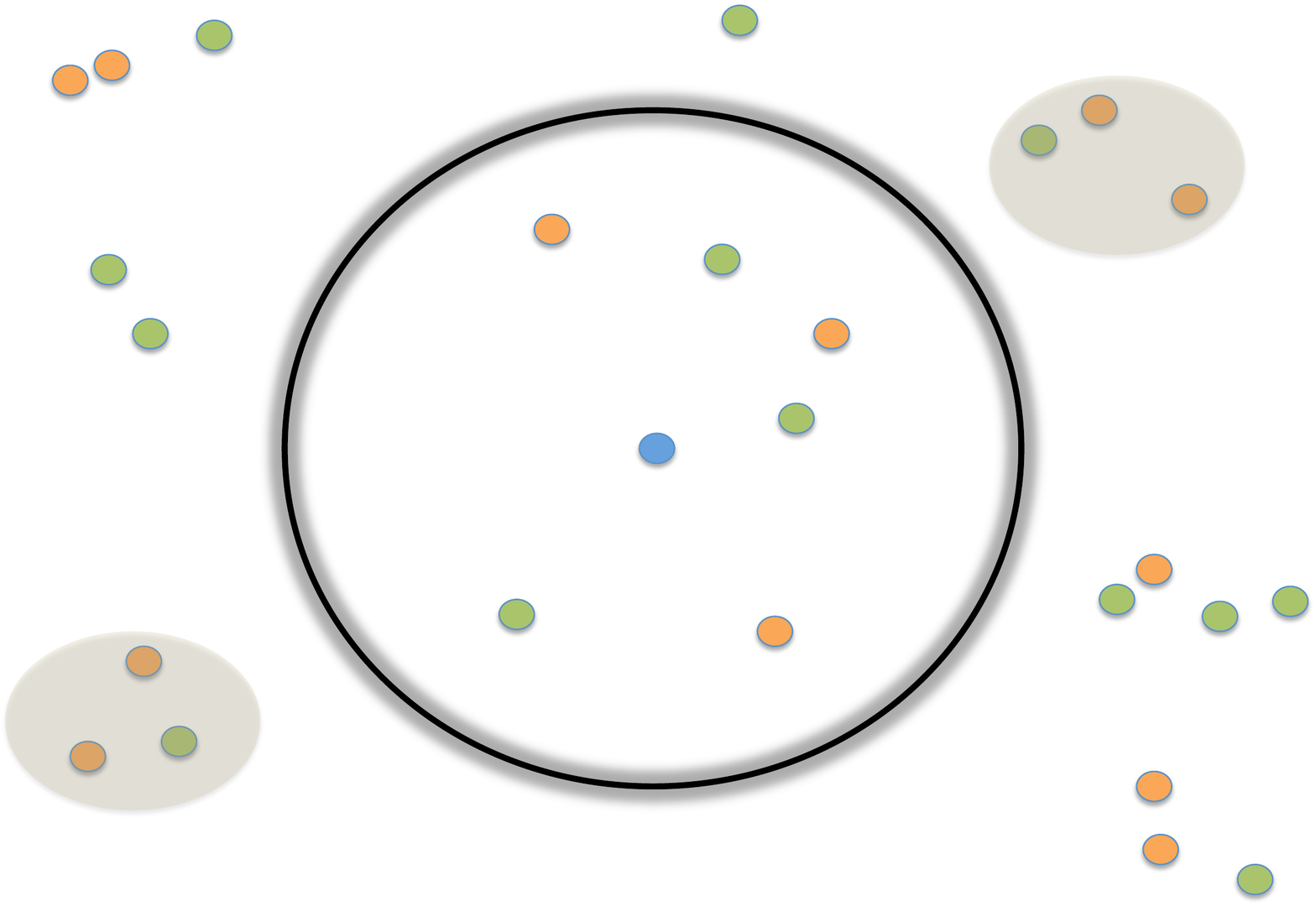}\\
\caption{Left: Helicoidal trajectories of ions and electrons submitted to a magnetic field. {\small (Source: S. Jardin's Lectures at Cemracs'10)}.
Right: Screen effect and Debye Length in a plasma.}
\label{figHtDl} 
\end{center}
\end{figure} 
As illustrated in left picture of Figure \ref{figHtDl}, because of the mass ratio between ions and electrons, the order of magnitude of electron Larmor Radius
is much smaller than that of ions.
In Tokamak, those order of magnitude are $\sim 5 \cdot 10^{-4} m$\ for electrons and $\sim  10^{-2} m$ for ions.
Those two lengths are two important scales of Tokamak Plasma Physics.

To simplify the purpose, from now we will forget the scale of electron Larmor Radius and only take ions under consideration.

\paragraph{Debye Length - } \label{tutu}We now explain an important phenomenon of Plasma Physics which is the Screen Effect, with its related length
which is  the Debye Length. 

It is well known that a charged particle generates around it an electric potential which depends on the inverse of the distance to it, say 
$\sim 1/r$.  Since this function is the Green Kernel of the opposite to the Laplace operator, when following the Mean Field routine
to deduce the model taking under consideration interactions of many particles, it is gotten that the electric field that fills the space is the solution to 
 $-\Delta \Phi=$ $\text{\it {Cst }{\tt Charge \!Density}}$.
 Because, of the long range of $\sim (1/r)$-potential, the regularization properties of minus Laplacian are such that space variations of the electric
 potential (and consequently of the electric field) occure at a large scale.
 
 Now, if the density of particles is relatively large, as it happens in Plasma, and if we consider a test particle (which is drawn in blue in right picture of
 Figure \ref{figHtDl}), the other particles (in the picture, orange particles are positively charged and green particles negatively charged) which are beyond 
 a given length $\lambda$, which is the Debye Length (the black circle in the picture is censed 
being of radius $\lambda$), may be gathered into subsets (the two grey ellipses in the picture are such subsets) which are such that the resulting action
of pairs of those subsets on the test particle is negligible.
This has the following consequence: it can be considered that the potential $\sim 1/r$ generated by a particle may be approximated by a  $\sim e^{-r/\lambda}/{r}$  
potential. Then,  following the Mean Field routine with this new assumption leads an electric field that fills the space which is the solution to 
 $-{\cal D}  \Phi = \text{\it {Cst }{\tt Charge \!Density}}$.
 Operator $-{\cal D}$ is the pseudo-differential operator which Green Kernel is $e^{-r/\lambda}/{r}$ whose regularization properties are much less than that of minus
 Laplacian.
 Then the regularity of the electric field is not as high as we could think before this formal analysis. In particular the electric field can be prone to variations over
 lengths ranging from several Debye Lengths to several hundreds of Debye Lengths.  Those variations are maybe not extremely large, but as no principle 
 fends them off, they are stable.
\paragraph{What Ions see -} In a very simplified slant, Electric field in a plasma is a bed with variations which are not
so large but not negligible, and, more than anything, stable at scales ranging from several Debye Lengths to several hundreds of Debye Lengths.\begin{figure}[htbp]
\begin{center}
\includegraphics[width=6cm, bb=1cm 5cm 20cm 25cm]{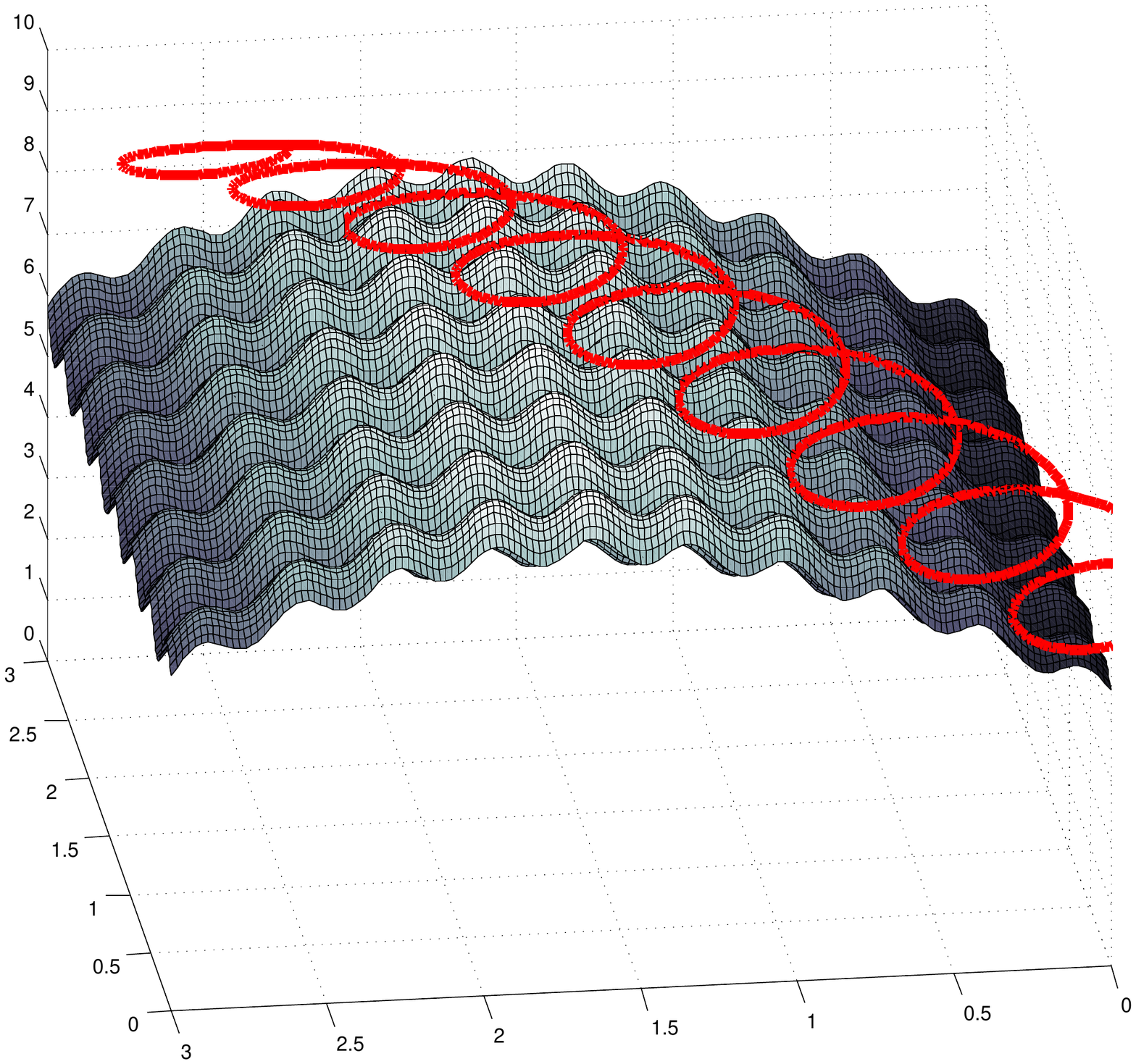}
\caption{Symbolic picture of oscillating trajectory of ion in a electric field with variations with typical size comparable
with trajectory oscillation.}
\label{figWis} 
\end{center}
\end{figure}
In Figure  \ref{figWis} we symbolically represent this bed as a charcoal grey surface. (It is symbolic because a surface is drawn while, in real world, Electric Field fills a 3D space.) 
This surface shows variations at two scales which are cenced to be about 100 times the Debye Lengths and 1.000 times the Debye Lengths.
Those variations are not of large amplitude but are stable.
Now, ion motion, symbolically represented in red in the figure, shows an oscillation of amplitude which is about 100 times the Debye Length
(ion Larmor Radius $\sim  10^{-2} m$, Debye Length : $\sim 10^{-4} m$).
Hence, because of the stability of bed variations, its resulting effect can be important on ion trajectory after several periods of oscillation.
To take into account this resulting effect, the concept of {gyro-average} was introduced. 
\paragraph{Dimensionless Dynamical System -}
Here, we do not present the scaling routine, we only mention that two small parameters may be introduced.
The first one is 
\begin{gather}\label{ddseq305} 
 \eps \sim \frac{\text{Ion Larmor Radius}}{\text{Tokamak size}} \sim \frac{10^{-2} m}{10 m}\sim 10^{-3}
\end{gather} 
The second one is $\eta$ the
Debye Length linked variations of the Electric Field and Potential. It is defined as
\begin{gather}\label{ddseq306} 
 \eta \sim \frac{\text{Characteristic scale of variations of the Electric Field}}{\text{Tokamak size}}.
\end{gather} 
If we follow a scaling routine, we will get that the dimensionless Electric Potential and Field write
\begin{gather}
\label{PhiE} 
\Phi(\xvec) = \Phi_0(\xvec) + \eta \Phi_1(\frac{\xvec}{\eta}) \text{ ~ and ~ }
   \Evec(\xvec)=\Evec_0(\xvec)+ \Evec_1(\frac{\xvec}{\eta}),
\end{gather} 
with
\begin{gather}
\label{Pot} 
\Evec(\xvec) = -(\nabla \Phi)(\xvec) ,  
\end{gather}
and that the dimensionless Magnetic Field write
\begin{gather}
\label{VecPot} 
\frac{\Bvec(\xvec)}{\eps} = (\nabla\times\frac{\Avec}{\eps})(\xvec).
\end{gather}
Then, the dimensionless trajectory, is solution to the following dynamical system:
         \begin{align} 
             &\fracp{\Xvec}{t}=\Vvec            
              \label{Dynsyst1} 
\\
             &\fracp{\Vvec}{t}=\Evec_0(\Xvec)+  \Evec_1(\frac{\Xvec}{\eta})+\Vvec \times \frac{\Bvec(\Xvec)}{\eps}.
             \label{Dynsyst2} 
          \end{align}
We guess that anybody agree with the fact that this system is a good model to describe motion of charged particle 
within Tokamak.
\paragraph{Gyro-kinetic Model -}Yet, there is the Gyro-Kinetic Model.  It claims that a trajectory $(\text{\Large $r$},\text{\Large $\psi$}, \text{\Large $\varphi$}, W_\parallel, J,\Gamma)$
in a coordinates system $(r,\psi, \varphi, w_\parallel, j, \gamma)$ we will discuss on later on, is solution to: 
\begin{gather}
\label{GyrMod} 
\begin{aligned}
&\ds \fracp{\text{\Large $r$}}{t}=(ED_r+ MCD_r), &
&\ds \fracp{\text{\Large $\psi$}}{t}=\frac{W_\parallel}{q(\text{\Large $r$})R}+ \frac{ED_\psi+ MCD_\psi}{\text{\Large $r$}},\\
&\ds \fracp{ \text{\Large $\varphi$}}{t}=\frac{W_\parallel}{R},&
&\ds \fracp{W_\parallel}{t}=(E_0+\langle\!E_1\!\rangle)_\parallel -\frac{J}{\eps} \Nabla_\parallel{|\Bvec|}+ \frac{W_\parallel}{|\Bvec|}ED\cdot(\nabla|\Bvec|),
\end{aligned}
\end{gather}
\begin{figure}[htbp]
\begin{center}
\includegraphics[width=11cm]{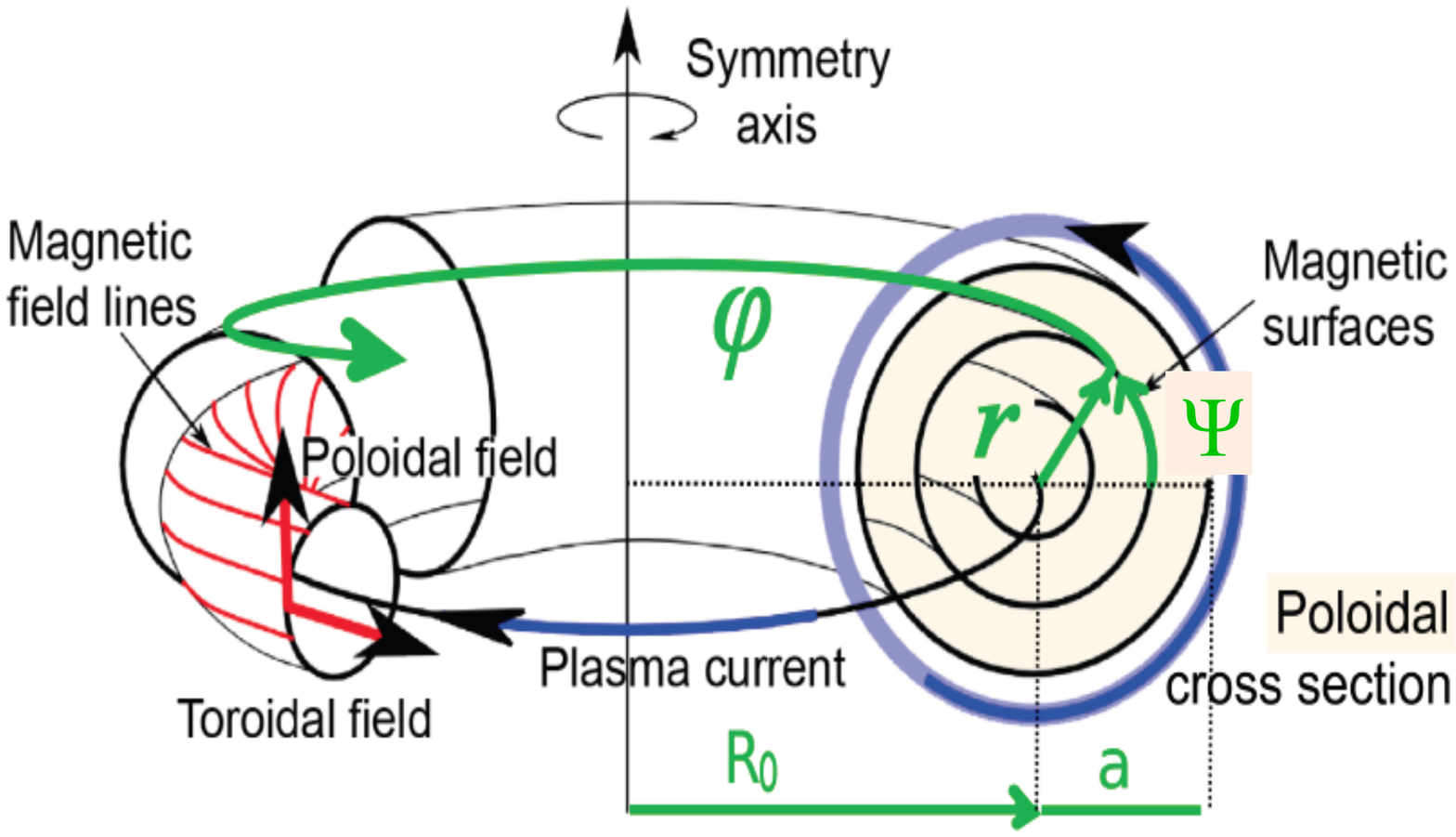}\vspace{-2.3cm}
\caption{Position coordinates in use in a Tokamak. {\small (Source: V. Grandgirard's Lectures at Cemracs'10).}}
\label{figTokamak} 
\vspace{-6mm}
\end{center}
\end{figure}
where
\begin{gather}
\label{GyrMod2} 
R=R_0+\text{\Large $r$}\cos(\text{\Large $\psi$}), ~ 
\ds ED= \eps\frac{\Bvec\times(E_0+\langle\!E_1\!\rangle)}{|\Bvec|^2}, ~ 
\ds MCD= \eps\frac{\Bvec\times(\nabla|\Bvec|)}{|\Bvec|^3}(W_\parallel^2+J\frac{|\Bvec|}{\eps}),
\end{gather}
and
\begin{gather}
\label{GyrMod3} 
\ds \Nabla_\parallel =  \frac{1}{R}\Big(\fracp{}{\varphi} + \frac{1}{q(\text{$r$})}\fracp{}{\psi} \Big).
\end{gather}

~

We now explain how to read equations \eqref{GyrMod} - \eqref{GyrMod3}.  A Tokamak is the interior of a torus with a vertical symmetry
axis. (A Tokamak and its associated position coordinate system  $(r,\psi, \varphi)$ is given in Figure \ref{figTokamak}.) Hence, it
is thought as  a small disc carried by  a large circle of radius $R_0$. Then, for a given point of the Tokamak, $r$ is (\emph{almost}) its distance from
the center of the small disk, $\psi$ is (\emph{almost}) the angle that the point makes in the small disk and $\varphi$ is (\emph{almost}) the angle that it makes on 
the large circle. 
\\
In what concerns the velocity variable system $(w_\parallel, j, \gamma)$,  if a moving point is in  $(r,\psi, \varphi)$, the $w_\parallel$-component 
of its velocity stands (\emph{almost}) for the projection of its velocity on the direction the Magnetic Field has in $(r,\psi, \varphi)$ (which is not far from being tangent to the large 
circle but which has a small component tangent to the small circle), $j$ is (\emph{almost}) half the square
of the norm of its velocity on the plan orthogonal to the Magnetic Field direction, and $\gamma$ is (\emph{almost}) the angle that this projection makes in this plan.

Above, precision "\emph{almost}" is essential, and constitutes the core of the Gyro-Kinetic Approximation which consists, as we will see, in building infinitesimal 
changes of coordinates to set out  a dynamical system shaped like \eqref{GyrMod}. 

We now comment system \eqref{GyrMod}, forgetting that the variables are only \emph{almost} what they are supposed to be. 
The first equation in \eqref{GyrMod}  says that the $r$-component of a particle trajectory varies with a velocity which is the projection on the direction of the vector
joining the center of the small disk and the point of the Electric Drift velocity $ED$ and the Magnetic Curvature Drift velocity $MCD$, both defined in \eqref{GyrMod2}.
The Electric Drift velocity $ED$ involves the cross product of the Magnetic Field and what is called the gyro-average of the Electric Field, which is indicated by symbol $\langle\;\rangle$. 
The Magnetic Curvature Drift velocity $MCD$ involves the cross product of the Magnetic Field with the gradient of the Magnetic Field norm.
\\
In the second equation of \eqref{GyrMod}, $q(r)$ is called the quality factor. Essentially, it is the number of revolutions in the small disk a particle makes while making
one revolution along the large circle. Then, term ${W_\parallel}/({q(\text{\Large $r$})R})$, where $R$ is defined in \eqref{GyrMod2},  generates the winding of the particle trajectory in the small disk (this is illustrated by the red
curves drawn on the left of the figure). The second term of the second equation of \eqref{GyrMod} is linked with the action of the Electric Drift velocity and Magnetic Curvature Drift velocity.
\\
The third equation involves the contribution ${W_\parallel}/{R}$ of the  $w_\parallel$-component of the particle velocity to the variation of the $\varphi$-component of the particle position.
\\
The last equation describes the variation of $w_\parallel$-component of the particle velocity. It involves, among other terms, the projection of the gyro-average of the Electric Field on the direction of the 
Magnetic Field and the variation of the Magnetic Field norm in the direction of the Magnetic Field,  materialized by operator $\Nabla_\parallel$ defined by  \eqref{GyrMod3}.

In system \eqref{GyrMod}, their is neither equation for $J$ (the $j$-component of the particle velocity) nor for $\Gamma$ (the $\gamma$-component of the particle velocity).
This may be considered as a strange fact. The reason why this happens is that the dynamics of $\text{\Large $r$},\text{\Large $\psi$}, \text{\Large $\varphi$}, W_\parallel$, 
the $r$-, $\psi$-,  $\varphi$-,  $w_\parallel$-components of the particle trajectory and velocity do not depend on the $\gamma$-component of the particle velocity and that
the its $j$-component is invariant. This is summarized as:
\begin{gather}
\label{GyrMod22} 
\begin{aligned}
&\ds \fracp{J}{t}=0, &
&\ds \fracp{\Gamma}{t}=\text{Does not matter}.
\end{aligned}
\end{gather}

\paragraph{What we explain and what we do not explain -} In this document we give the key ideas that bring the deduction of the Gyro-Kinetic 
Dynamical System \eqref{GyrMod2}
from the original one \eqref{Dynsyst1} - \eqref{Dynsyst2}.

Among thing that we do not tackle, there is the question of the specific position coordinates $(r,\psi, \varphi)$. We will remain with a generic  position coordinate system.
We will not talk about Gyro-average either. 
The question of the coupling of Gyro-Kinetic Approximation with computation of Electric and Magnetic Fields, using Quasi-Neutral Poisson Equation or Maxwell Equation
is also not addressed here.

To end this paragraph we mention that the method only works in case when $\eta$ is larger than $\eps$, or in other words,
when $\eta= \eps^{1-\kappa}$ for $\kappa>0$. This is a important limitation in the theory.
\label{gghh}
%
%
\section{Methode summarize}
\paragraph{Key result -} The method allowing us  to go from system \eqref{Dynsyst1} - \eqref{Dynsyst2} 
to a system of the kind of \eqref{GyrMod},
i.e. to make disappear two equations, has nothing to do with magic. It is based on the following result:
\begin{theorem}
\label{KrTh} 
If, in  a given coordinate system $\rvec =(r_1, r_2, r_3, r_4, r_5, r_6),$ a Hamiltonian Dynamical System writes:
\begin{gather} 
\label{HamSystSh1} 
              \ds\fracp{\Rvec}{t} = \Pcal(\Rvec)\Nabla_{\rvec} H(\Rvec), ~~~~
              \Pcal(\Rvec) = \left(\begin{array}{c|cc}
                    \text{\huge${\cal M}$} & \begin{array}{c}0 \vspace{-3pt}\\ : \\0 \end{array} & \begin{array}{c}0 \vspace{-3pt}\\ : \\0 \end{array} \\
                    \hline
                     0 \cdots 0& 0 & 1\\
                     0 \cdots 0& -1 & 0 \\
                                   \end{array}\right),
\end{gather}
(for a trajectory $\Rvec =(R_1, R_2, R_3, R_4, R_5, R_6)^T$))
with a Hamiltonian Function that does not depend on the last variable, i.e.
\begin{gather}
\label{HamSystSh2} 
\fracp{H}{r_6} =0,
\end{gather}
Then, submatrix ${\cal M}$ does not depend on the two last variables, i.e. 
\begin{gather}
\label{HamSystSh3} 
   \ds \fracp{{\cal M}}{r_5}=0 \text{ and } \fracp{{\cal M}}{r_6}=0,
\end{gather}
and consequently,  time-evolution of the four first components $R_1, R_2, R_3, R_4$ is independent
of the last component $R_6$;
and, the penultimate component $R_5$ of the trajectory in not time-evolving, i.e.
\begin{gather}
\label{HamSystSh4} 
  \ds \ds \fracp{R_5}{t}=0.
\end{gather}
\end{theorem}
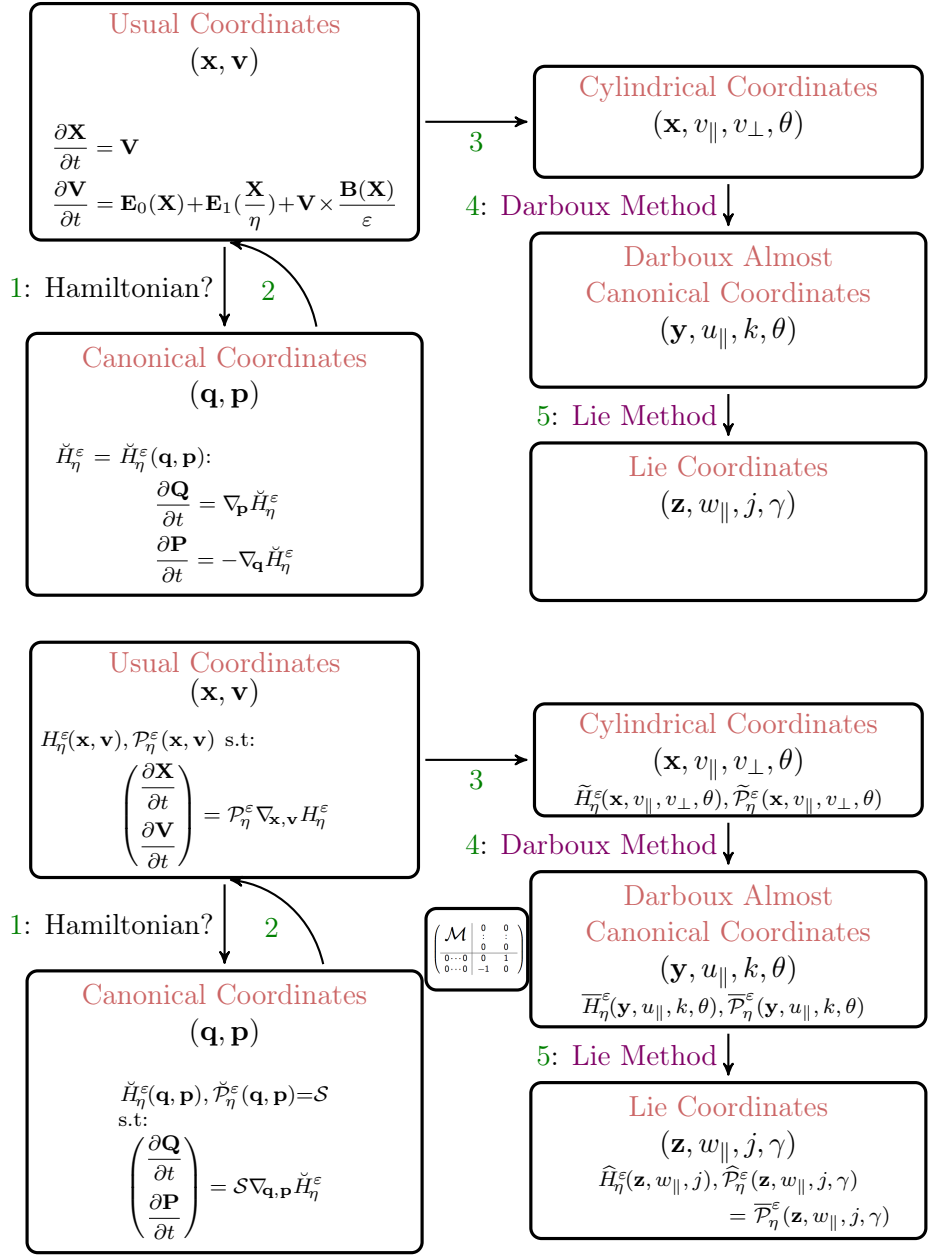
\begin{figure}[htbp]
\begin{center}
\begin{tikzpicture}[node distance=1cm, auto,]
  \node[punkt] (coordusuelles) {{\color{light-red}Usual Coordinates} \\
          $(\xvec,\vvec)$ \\ ~ \vspace{-2mm} \\
        \scriptsize $~$~~~~~~~~~~~~~~~~~~~\vspace{-3mm}
          \scriptsize \begin{align*} 
             &\fracp{\Xvec}{t}=\Vvec\\
             &\fracp{\Vvec}{t}=\Evec_0(\Xvec)\!+\! \Evec_1(\frac{\Xvec}{\eta})\!+\!\Vvec \!\times\! \frac{\Bvec(\Xvec)}{\eps}
          \end{align*}
  };
   \node[punkt, inner sep=5pt,below=1.2cm of coordusuelles] (coordcanon) {{\color{light-red}Canonical Coordinates} \\ 
        $(\qvec,\pvec)$ \\ ~ \\
        \scriptsize $\breve H^\eps_{\!\eta}=\breve H^\eps_{\!\eta}(\qvec,\pvec)$:~~~~~~~~~~~~~~~~~~~\vspace{-3mm}
         \scriptsize\begin{align*} 
             &\fracp{\Qvec}{t}=\Nabla_\pvec \breve H^\eps_{\!\eta} \\ 
             &\fracp{\Pvec}{t}=-\Nabla_\qvec \breve H^\eps_{\!\eta}
         \end{align*}
    }
    edge[pil,bend right=35] node[auto] {{\color{vert}2}}(coordusuelles.south)
    edge[pil,<-,] node[auto] {{\color{vert}1}:  Hamiltonian?}(coordusuelles.south);
     \node[punkt,right=1.5cm of coordusuelles] (coordcyl) {{\color{light-red} Cylindrical Coordinates} \\ 
            $(\xvec,v_\parallel,v_\perp,\theta)$\\
            \scriptsize$~~~$
     }
     edge[pil,<-,] node[auto] {{\color{vert}3}}(coordusuelles.east);
      \node[punkt,inner sep=5pt,below=0.7cm of coordcyl] (coorddarboux) {{\color{light-red} Darboux Almost Canonical Coordinates} \\
            $(\yvec,u_\parallel,k,\theta)$\\
            \scriptsize$~~~$
     }
     edge[pil,<-,] node[auto] {{\color{vert}4}: {\color{violet}Darboux Method}}(coordcyl.south);
     \node[punkt,inner sep=5pt,below=0.7cm of coorddarboux] (coorddarlie) {{\color{light-red} Lie  Coordinates} \\
            $(\zvec,w_\parallel,j,\gamma)$\\
            \scriptsize$~~~$\\
            \scriptsize$~~~$\\
            \scriptsize$~~~$

     }
     edge[pil,<-,] node[auto] {{\color{vert}5}: {\color{violet}Lie Method}}(coorddarboux.south);
\end{tikzpicture}
\\~\\
\begin{tikzpicture}[node distance=1cm, auto,]
  \node[punkt] (coordusuelles) {{\color{light-red}Usual Coordinates}\vspace{-1mm} \\
          $(\xvec,\vvec)$ \\ ~ \vspace{-2mm} \\
         \scriptsize $H^\eps_{\!\eta}\!(\xvec,\vvec),\Pcal^\eps_{\!\eta}(\xvec,\vvec)$~s.t:~~~~~~~~~~~~~~~~~~~\vspace{-3mm}
         \scriptsize\begin{gather*} 
             \begin{pmatrix} \ds\fracp{\Xvec}{t} \\ ~\vspace{-1.5mm}\\ \ds \fracp{\Vvec}{t} \end{pmatrix} =
             \Pcal^\eps_{\!\eta}\Nabla_{\xvec,\vvec} H^\eps_{\!\eta}
         \end{gather*}
  };
   \node[punkt, inner sep=5pt,below=1.2cm of coordusuelles] (coordcanon) {{\color{light-red}Canonical Coordinates} \\ 
        $(\qvec,\pvec)$ \\ ~ \\
                \scriptsize $\breve H^\eps_{\!\eta}\!(\qvec,\pvec),\breve\Pcal^\eps_{\!\eta}(\qvec,\pvec)\text{$=$}\Scal$
                ~s.t:~~~~~~~~~~~~~~~~~~~\vspace{-3mm}
         \scriptsize\begin{gather*} 
             \begin{pmatrix} \ds\fracp{\Qvec}{t} \\ ~\vspace{-1.5mm}\\ \ds \fracp{\Pvec}{t} \end{pmatrix} =
             \Scal \Nabla_{\qvec,\pvec} \breve H^\eps_{\!\eta}
         \end{gather*}
    }
    edge[pil,bend right=35] node[auto] {{\color{vert}2}}(coordusuelles.south)
    edge[pil,<-,] node[auto] {{\color{vert}1}:  Hamiltonian?}(coordusuelles.south);
     \node[punkt,right=1.5cm of coordusuelles] (coordcyl) {{\color{light-red} Cylindrical Coordinates}\\ 
            $(\xvec,v_\parallel,v_\perp,\theta)$\\
            \scriptsize$\widetilde H^\eps_{\!\eta}\!(\xvec,v_\parallel,v_\perp,\theta),
            \widetilde\Pcal^\eps_{\!\eta}(\xvec,v_\parallel,v_\perp,\theta)$\vspace{-1mm}
     }
     edge[pil,<-,] node[auto] {{\color{vert}3}}(coordusuelles.east);
      \node[punkt,inner sep=5pt,below=0.7cm of coordcyl] (coorddarboux) {{\color{light-red} Darboux Almost Canonical Coordinates} \\
            $(\yvec,u_\parallel,k,\theta)$\\
            \scriptsize$\overline H^\eps_{\!\eta}\!(\yvec,u_\parallel,k,\theta),\overline \Pcal^\eps_{\!\eta}(\yvec,u_\parallel,k,\theta)$ \vspace{-1mm}
     }
     edge[pil,<-,] node[auto] {{\color{vert}4}: {\color{violet}Darboux Method}}(coordcyl.south);
     \node[punkt,inner sep=5pt,below=0.7cm of coorddarboux] (coorddarlie) {{\color{light-red} Lie  Coordinates} \\
            $(\zvec,w_\parallel,j,\gamma)$\\
            \scriptsize$\widehat H^\eps_{\!\eta}\!(\zvec,w_\parallel,j),\widehat \Pcal^\eps_{\!\eta}(\zvec,w_\parallel,j,\gamma)$\\
            ~~~~~~~~~~~~~~~~~~$=\overline \Pcal^\eps_{\!\eta}(\zvec,w_\parallel,j,\gamma)$    \vspace{-0.7mm}
     }
     edge[pil,<-,] node[auto] {{\color{vert}5}: {\color{violet}Lie Method}}(coorddarboux.south);
     \node[punktt,inner sep=5pt,left=0.0mm of coorddarboux] (matrixshape) {~\hspace{-2mm}\includegraphics[width=12mm]{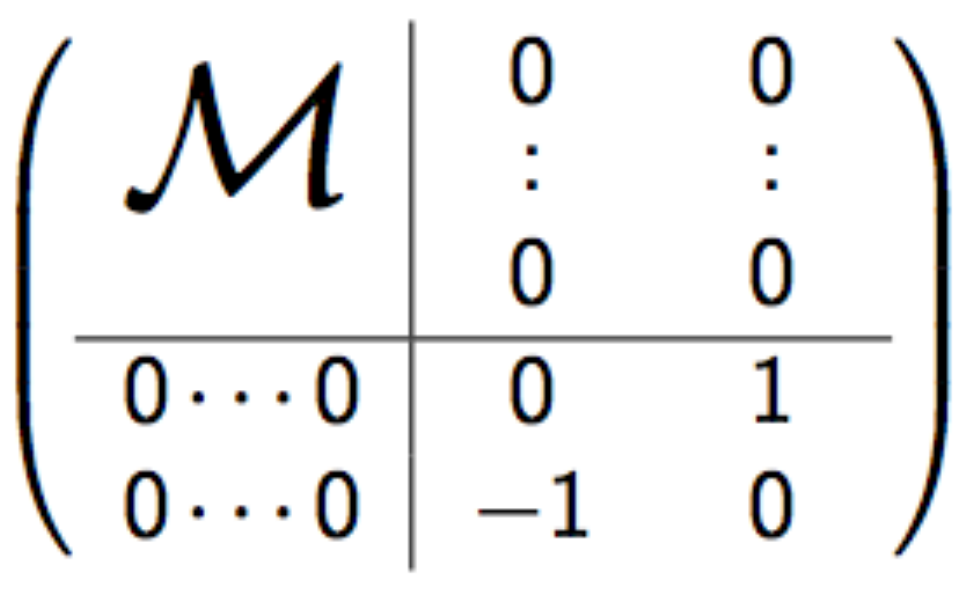} };
\end{tikzpicture}
\caption{The method is made of 5 steps. 1: Check that the Dynamical System \eqref{Dynsyst1}- \eqref{Dynsyst2} is well
     Hamiltonian. 2: Write the system using Hamiltonian Function and Poisson Matrix. 3: Write the system in a cylindrical in velocity coordinate system 
     using the formula giving how the Poisson Matrix and the Hamiltonian Function is transformed by change of coordinates. 4: Make another change
     of coordinates in order to have the Poisson Matrix form allowing the application of the Key Result (Theorem \ref{KrTh}). 5: Make a last change of
     coordinates, leaving the Poisson Matrix form unchanged and leading to a Hamiltonian independent of the last variable.}
\label{figPanorama} 
\end{center}
\end{figure}

The proof of property \eqref{HamSystSh3} is relatively straightforward by a calculation.
Equation  \eqref{HamSystSh4} is a direct consequence of  \eqref{HamSystSh2} that yields a sixth component
of $\Nabla_{\rvec} H(\Rvec)$ with worth 0 and consequently, because of the form of  $\Pcal(\Rvec)$ a fifth line of Dynmical
System \eqref{HamSystSh1} which is exactly  \eqref{HamSystSh4}.

Essentially, the Gyro-Kinetic version of Dynamical System \eqref{Dynsyst1} - \eqref{Dynsyst2} exhibit essentially a property like
\eqref{HamSystSh4} and the property that  the four first components of the trajectory are independent of the last component.

~

As a matter of fact, building the Gyro-Kinetic Approximation consists in building a change of coordinates that bring a system of
coordinates in which $\Pcal(\Rvec)$ has the same form as in  \eqref{HamSystSh1} and in which \eqref{HamSystSh2} is true.
\paragraph{Panorama -} As illustrated in Figure  \ref{figPanorama},  the method to build the desired change of coordinates is made of 5 steps. 
The first one consists in checking that the Dynamical System  \eqref{Dynsyst1} - \eqref{Dynsyst2} is well Hamiltonian. This is symbolized by arrow 1 in the
top picture of the Figure.  Once this is done, we can go back into the Usual Coordinate System but knowing that the system writes as in the
square which is in the top-left  of the bottom picture, i.e. involving a Poisson Matrix $\Pcal^\eps_{\!\eta}$ and the gradient $\Nabla_{\xvec,\vvec} H^\eps_{\!\eta}$ of a Hamiltonian
Function.  It may be written in that form in any Coordinate System, and formula give how to transform
the Hamiltonian Function and the Poisson Matrix while changing of Coordinates.
The goal of the third step is to introduce a Cylindrical in Velocity Coordinate system which is known has being close to the Gyro-Kinetic Coordinate System. Writing the system in this
coordinate system uses the formula giving how the Poisson Matrix and the Hamiltonian Function are transformed by change of coordinates.
(At this level, it would be also possible to make a position change of coordinates fitting Tokamak geometry.)
 In the fourth step, we make another change of coordinates in order to have the Poisson Matrix form allowing the application of the Key Result (Theorem \ref{KrTh}). 
 In the fifth step, we make a last change of coordinates, leaving the Poisson Matrix form unchanged and leading to a Hamiltonian independent of the last variable.
%
%
\section{Hamiltonian System}
We make here the first step: we check that Dynamical System  \eqref{Dynsyst1} - \eqref{Dynsyst2} is Hamiltonian. Essentially, this means nothing but
that there exists a coordinate system in which it writes:

\begin{align}
              \label{HHSS676} 
             &\fracp{\Qvec}{t}=\Nabla_\pvec \breve H^\eps_{\!\eta}, \\ 
              \label{HHSS677} 
             &\fracp{\Pvec}{t}=-\Nabla_\qvec \breve H^\eps_{\!\eta},
\end{align}
or
\begin{gather} 
\label{HHSS678} 
             \begin{pmatrix} \ds\fracp{\Qvec}{t} \\ ~\vspace{-1.5mm}\\ \ds \fracp{\Pvec}{t} \end{pmatrix} =
             \Scal \Nabla_{\qvec,\pvec} \breve H^\eps_{\!\eta},
             \text{ ~ with ~ } \Scal = \begin{pmatrix} 0 & I_3  \\ -I_3 & 0 \end{pmatrix},
\end{gather}
where $\breve H^\eps_{\!\eta} $ is a function.

\paragraph{The Canonical Coordinates -}
If we use \eqref{Pot} and \eqref{VecPot} in \eqref{Dynsyst1} - \eqref{Dynsyst2}, we get that in the
Usual Coordinates $(\xvec,\vvec) = (x_1,x_2,x_3, v_1,v_2,v_3)$, 
trajectory : $(\Xvec(t;\xvec,\vvec,s), \Vvec(t;\xvec,\vvec,s))$
is solution to
\begin{align}
\label{DynSystPot1} 
&\ds \fracp{\Xvec}{t}=\Vvec,  \\
\label{DynSystPot2} 
&\ds \fracp{\Vvec}{t}=  - \nabla \big[\Phi_0(\Xvec) + {\eta} \Phi_1(\frac{\Xvec}{\eta}) \big] +\Vvec \times \frac{\nabla\times\Avec(\Xvec)}{\eps}.
\end{align}
We recall that {$(\Xvec,\Vvec) = (X_1,X_2,X_3, V_1,V_2,V_3)$}.\\

We will show, in the next paragraph, that in coordinate system $(\qvec,\pvec) = (q_1,q_2,q_3, p_1,p_2,p_3)$ defined by 
\begin{gather}
\label{VerifHamilt101} 
\qvec=\xvec,   ~~ \ds\pvec=\vvec+\frac{\Avec(\xvec)}{\eps},
\end{gather}
with reverse transformation given by
\begin{gather}
\label{VerifHamilt102} 
\xvec=\qvec, ~~ \ds\vvec=\pvec-\frac{\Avec(\qvec)}{\eps},
\end{gather}
the trajectory $(\Qvec(t;\qvec,\pvec,s), \Pvec(t;\qvec,\pvec,s))$ ($(\Qvec,\Pvec) = (Q_1,Q_2,Q_3, P_1,P_2,P_3)$) which is 
given by
\begin{gather}
\label{VerifHamilt103} 
\Qvec=\Xvec, ~~ \ds\Pvec=\Vvec+\frac{\Avec(\Xvec)}{\eps}, ~~~~~~\Xvec=\Qvec, ~~\ds \Vvec=\Pvec-\frac{\Avec(\Qvec)}{\eps},
\end{gather}
is solution of a Dynamical System of the form  \eqref{HHSS676}  -  \eqref{HHSS677} or  \eqref{HHSS678},  with 
\begin{gather}
\label{VerifHamilt104} 
\ds\breve H^\eps_{\!\eta}(\qvec,\pvec) = \frac 12 \Big| \pvec - \frac{\Avec(\qvec)}{\eps} \Big|^2 + \Phi_0(\qvec) + {\eta} \Phi_1(\frac{\qvec}{\eta}).
\end{gather}
Those $(\qvec,\pvec)$ coordinates will be called Canonical Coordinates.

\paragraph{Check of Canonical nature of Canonical Coordinates -}
Making calculations suggested by \eqref{HHSS676}  -  \eqref{HHSS677} with $\breve H^\eps_{\!\eta}$  given by  \eqref{VerifHamilt104}
gives
\begin{align}
\label{ChCanNat301} 
&\ds \fracp{\Qvec}{t}=\Nabla_\pvec \breve H^\eps_{\!\eta}(\Qvec,\Pvec)=\Pvec - \frac{\Avec(\Qvec)}{\eps}, \\
\label{ChCanNat302} 
&\ds \fracp{\Pvec}{t}=-\Nabla_\qvec \breve H^\eps_{\!\eta}(\Qvec,\Pvec)=
     \frac{(\nabla\Avec(\Qvec))^{\!T}}{\eps}\Big(\Pvec-\frac{\Avec(\Qvec)}{\eps}\Big) 
    - \nabla \big[\Phi_0(\Qvec) + {\eta} \Phi_1(\frac{\Qvec}{\eta}) \big],
\end{align}
which applying the following formula
\begin{gather}
\label{ChCanNat303} 
  (\nabla\Avec)^{\!T}(\pvec-\Avec) =(\nabla\Avec)(\pvec-\Avec) + (\pvec-\Avec) \times (\nabla\times\Avec),
\end{gather}
yields
\begin{align}
\label{ChCanNat304} 
&\ds \fracp{\Qvec}{t}=\Pvec - \frac{\Avec(\Qvec)}{\eps},\\
\label{ChCanNat305} 
&\ds \fracp{\Pvec}{t}- \frac{(\nabla\Avec(\Qvec))}{\eps}\Big(\Pvec-\frac{\Avec(\Qvec)}{\eps}\Big) =
 \Big(\Pvec-\frac{\Avec(\Qvec)}{\eps}\Big) \times \frac{\nabla\times\Avec(\Qvec)}{\eps}
   - \nabla \big[\Phi_0(\Qvec) + {\eta} \Phi_1(\frac{\Qvec}{\eta}) \big].
\end{align}
Now,using \eqref{VerifHamilt103} we get
 \begin{align}
\label{ChCanNat307} 
&\ds \fracp{\Xvec}{t}=\Vvec ,\\
&\ds~\hspace{2.2cm}~ \fracp{\Pvec}{t}- \frac{(\nabla\Avec(\Qvec))}{\eps}\Big(\fracp{\Qvec}{t} \Big) =\fracp{\ds\Big[\Pvec-\frac{\Avec(\Qvec)}{\eps}\Big]}{t}= 
\nonumber \\
\label{ChCanNat308} 
&\fracp{\Vvec}{t}=
\Big(\Pvec-\frac{\Avec(\Qvec)}{\eps}\Big) \times \frac{\nabla\times\Avec(\Qvec)}{\eps}
   - \nabla \big[\Phi_0(\Qvec) + {\eta} \Phi_1(\frac{\Qvec}{\eta}) \big]\\
\nonumber
& ~ \hspace{5mm}  =\Vvec \times \frac{\nabla\times\Avec(\Xvec)}{\eps} - \nabla \big[\Phi_0(\Xvec) + {\eta} \Phi_1(\frac{\Xvec}{\eta}) \big] . 
\end{align}
which is \eqref{DynSystPot1} -  \eqref{DynSystPot2}, and then proving the equivalence between the two systems and leads to the conclusion
that the Dynamical System we take interest in is Hamiltonian. 

%

\paragraph{As by-products : Poisson Matrix, Poisson Bracket, Change of Coordinates Formula -}
To end the first step we now list some properties we can moreover deduce from the fact that Dynamical System  \eqref{Dynsyst1} - \eqref{Dynsyst2} 
is Hamiltonian.\\

In any coordinate system $\rvec =(r_1, r_2, r_3, r_4, r_5, r_6),$ the Dynamical System writes:
\begin{gather} 
\label{PmPbCcf24} 
              \ds\fracp{\Rvec}{t} = \Pcal(\Rvec)\Nabla_{\rvec} H(\Rvec),
\end{gather}
for a Poisson Matrix $\Pcal$, which is antisymmetric, and a Hamiltonian Function $H$.\\ 

We can consider the Poisson Bracket which is defined for two regular functions $f$ and $g: \rit^6 \rightarrow \rit$ by
\begin{gather} 
\label{PmPbCcf25} 
\ds \poibrack{f}{g}(\rvec) = (\Nabla_{\rvec}f(\rvec))\cdot (\Pcal(\rvec)(\Nabla_{\rvec}g(\rvec))).
\end{gather}
If one of the two function is vector-valued, i.e. if for instance $\fvec: \rit^6 \rightarrow \rit^6$ and $g:\rit^6 \rightarrow \rit$,
then it is defined, for any component $i= 1,\dots, 6$, by
\begin{gather} 
\label{PmPbCcf26} 
\ds (\poibrack{\fvec}{g}(\rvec))_i = (\Nabla_{\rvec}\fvec_i(\rvec))\cdot (\Pcal(\rvec)(\Nabla_{\rvec}g(\rvec))).
\end{gather}
Then, if we introduce the Coordinate Function $\Ivec$, defined by
$\ds\Ivec(\rvec) = \rvec$,
then \eqref{PmPbCcf24} reads also:
\begin{gather} 
\label{PmPbCcf009} 
              \ds\fracp{\Rvec}{t}  = \poibrack{\Ivec}{H}(\Rvec).
\end{gather}
\begin{remark}
In papers of physicists, formula  \eqref{PmPbCcf009} is generally read $\ds   \ds\fracp{\rvec}{t}  = \poibrack{\rvec}{H}$.
\end{remark}

In another coordinate system $\tilde\rvec =(\tilde r_1, \tilde r_2, \tilde r_3, \tilde r_4, \tilde r_5, \tilde r_6)$ 
which is tied with the first one by change-of-coordinates formula $\tilde\rvec = \rtrf(\rvec)$ and $\rvec=\tilde\rtrf(\tilde\rvec)=\rtrf^{-1}(\tilde\rvec)$,
the system writes
\begin{gather} 
\label{PmPbCcf27} 
              \ds\fracp{\tilde\Rvec}{t} = \tilde\Pcal(\tilde\Rvec)\Nabla_{\tilde\rvec} \tilde H(\tilde\Rvec),
\end{gather}
where $\tilde\Rvec = \rtrf(\Rvec)$, with Hamiltonian Function and Poisson Matrix given by
\begin{gather}
\label{PmPbCcf28} 
\ds \tilde H(\tilde\rvec) =  H(\tilde\rtrf(\tilde\rvec)) \text{ ~ and ~ }
\ds (\tilde\Pcal(\tilde\rvec))_{ij}=\poibrack{\rtrf_i}{\rtrf_j}(\tilde\rtrf(\tilde\rvec)).
\end{gather}

\begin{remark}
From \eqref{VerifHamilt104} and  applying  \eqref{PmPbCcf28}, in Usual Coordinates, Hamiltonian Function  expression is 
\begin{gather}
\label{HfPmUc42} 
\ds H^\eps_{\!\eta}\!(\xvec,\vvec)= \frac 12 \big| \vvec \big|^2 + \Phi_0(\xvec) + {\eta} \Phi_1(\frac{\xvec}{\eta}),
\end{gather}
leading to 
\begin{gather}
\label{HfPmUc43} 
\ds \Nabla_{\xvec,\vvec} H^\eps_{\!\eta}\! = \begin{pmatrix} \nabla\big[ \Phi_0(\xvec) + {\eta} \Phi_1(\frac{\xvec}{\eta})\big] \\  \vvec \end{pmatrix},
\end{gather}
and Poisson Matrix expression is
\begin{gather}
\label{HfPmUc45} 
\ds \Pcal^\eps_{\!\eta}(\xvec,\vvec) =
\begin{pmatrix}
0 & I_3 \\
-I_3 & \frac{(\nabla A(\xvec))^T - (\nabla A(\xvec))}{\eps}
\end{pmatrix}.
\end{gather}
This can be easily checked calculating
\begin{gather}
\label{HfPmUc46} 
\Pcal^\eps_{\!\eta}(\xvec,\vvec)  \Nabla_{\xvec,\vvec} H^\eps_{\!\eta}\! = \begin{pmatrix}
\vvec \\
\ds -\nabla\big[ \Phi_0(\xvec) + {\eta} \Phi_1(\frac{\xvec}{\eta})\big]   + \frac{(\nabla A(\xvec))^T - (\nabla A(\xvec))}{\eps} \vvec
\end{pmatrix},
\end{gather}
using \eqref{ChCanNat302},and comparing with  \eqref{DynSystPot1} -  \eqref{DynSystPot2}.

Expression \eqref{HfPmUc42} is consistent with  \eqref{PmPbCcf28}, let us check that it is the same for \eqref{HfPmUc45}.
To do this, let us write the change of coordinates in the following way:
$(\xvec,\vvec)=\xtrfb(\qvec,\pvec)$ and $(\qvec,\pvec)=\ptrf(\xvec,\vvec)$.  Then, for instance: 
$\xtrfb_5(\qvec,\pvec) = p_2-\frac{\Avec_2(\qvec)}{\eps}$,  $\xtrfb_6(\qvec,\pvec) = p_3-\frac{\Avec_3(\qvec)}{\eps}$ and
\begin{gather}
\label{HfPmUc47} 
\nabla\xtrfb_5(\qvec,\pvec)  =
\begin{pmatrix} \frac1\eps \fracp{\Avec_2}{q_1}(\qvec) \\ \frac1\eps\fracp{\Avec_2}{q_2}(\qvec) 
           \\ \frac1\eps\fracp{\Avec_2}{q_3}(\qvec) \\ 0 \\ 1 \\ 0 \end{pmatrix}, ~ 
\ds \Scal = \begin{pmatrix} 0 & I_3  \\ -I_3 & 0 \end{pmatrix}, ~ 
\nabla\xtrfb_6(\qvec,\pvec)  =
\begin{pmatrix} \frac1\eps \fracp{\Avec_3}{q_1}(\qvec) \\ \frac1\eps\fracp{\Avec_3}{q_2}(\qvec) 
            \\ \frac1\eps\fracp{\Avec_3}{q_3}(\qvec) \\ 0 \\ 0 \\ 1 \end{pmatrix}.
 \end{gather}
 Hence, since $\poibrack{\xtrfb_5}{\xtrfb_6} =  (\nabla\xtrfb_5)\cdot (\Scal\nabla\xtrfb_6)$, expression of  $(\Pcal^\eps_{\!\eta})^{~}_{56}$,
 obtained applying \eqref{PmPbCcf28}, is
 \begin{gather}
\label{HfPmUc48} 
\ds \poibrack{\xtrfb_5}{\xtrfb_6}(\ptrf(\xvec,\vvec)) = \frac1\eps \Big( \fracp{\Avec_3}{q_2}(\xvec) - \fracp{\Avec_2}{q_3}(\xvec) \Big),
\end{gather}
which is also the expression of  $(\Pcal^\eps_{\!\eta})^{~}_{56}$ obtained applying  \eqref{HfPmUc45}. 
(This may of course be led for the other entries of the matrix also.)
\end{remark}
%
%
%
\section{Cylindrical Coordinates}
\paragraph{Cylindrical Coordinates in velocity -}
We now turn to the second step which consists in setting the expression of Dynamical System  \eqref{DynSystPot1} - \eqref{DynSystPot2}
(or \eqref{Dynsyst1} - \eqref{Dynsyst2} or  \eqref{HHSS676}, \eqref{HHSS677},  \eqref{VerifHamilt104}) in a Cylindrical in velocity Coordinate
System $(\xvec,v_\parallel,v_\perp,\theta)$ which is such that
\begin{figure}[htbp]
\begin{center}
\includegraphics[width=11cm, bb=0 120 842 550]{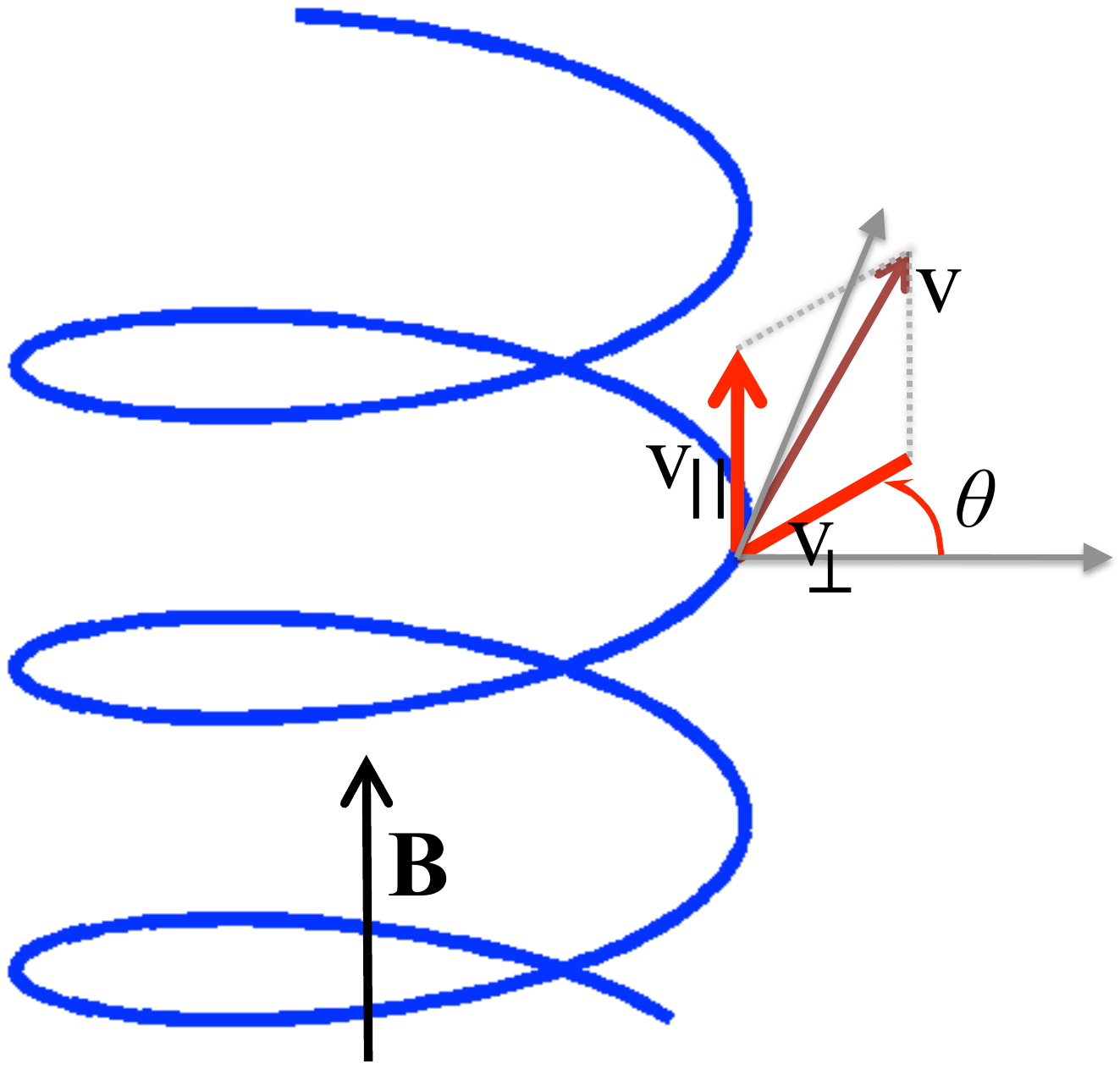} \vspace{-3mm} 
\caption{Helicoidal trajectory and  Cylindrical Coordinates for the velocity variable.}
\label{fig1TT} 
\end{center}
\end{figure}
\begin{gather}
\label{CyCooIVel789} 
\ds v_\parallel = \vvec\cdot \frac{\Bvec}{|\Bvec|}, ~
\ds v_\perp =\Big|\vvec- \Big(\vvec\cdot \frac{\Bvec}{|\Bvec|}\Big)\frac{\Bvec}{|\Bvec|}\Big|, ~
\theta \text{ s.t. } \ds \vvec- \Big(\vvec\cdot \frac{\Bvec}{|\Bvec|}\Big)\frac{\Bvec}{|\Bvec|} = v_\perp (\cos\theta,\sin\theta).
\end{gather}
This means, as shown in Figure \ref{fig1TT} where a helicoidal trajectory induced by a Magnetic Field pointing upward,
that $v_\parallel$ is the projection of the velocity on the direction of the Magnetic Field, $v_\perp$ is the norm of the projection of the 
velocity on the plan orthogonal to the Magnetic Field and $\theta$ is the angle that this projection makes in this plan.
%

\paragraph{Hamiltonian Function and Poisson Matrix in Cylindrical Coordinates~-}
Applying formula \eqref{HfPmUc42}, we get the expression of the Hamiltonian Function in this System:
\begin{gather}
\label{CyCooIVel790} 
\ds \widetilde H^\eps_{\!\eta}\!(\xvec,v_\parallel,v_\perp,\theta)= 
\frac 12 \big( v_\parallel ^2+v_\perp^2\big) + \Phi_0(\xvec) + {\eta} \Phi_1(\frac{\xvec}{\eta}).
\end{gather}
We can also get the expression of the Poisson Matrix, which is heavy, applying formula \eqref{HfPmUc42}.\\

Here, we give its expression in the case when 
\begin{gather}
\label{CyCooIVel791} 
\Bvec(\xvec) = \begin{pmatrix} b(\xvec) \\ 0 \\ 0 \end{pmatrix}, ~  b>0,
\end{gather}
only. It is:
\begin{gather}
\label{CyCooIVel792} 
\widetilde \Pcal^\eps_{\!\eta}(\xvec,v_\parallel,v_\perp,\theta)=
\left(\begin{array}{cccccc}
0 & 0 & 0  & \frac{b(\xvec)}{\eps} & 0 & 0 \\
0 & 0 & 0 & 0 & -\sin(\theta) & -\frac{\cos(\theta)}{v_\perp} \\
0 & 0 & 0 & 0 & -\cos(\theta)& \frac{\sin(\theta)}{v_\perp} \\
 -\frac{b(\xvec)}{\eps} & 0 & 0 & 0  & \$\$ & \$ \\
0 & \sin(\theta) & \cos(\theta) & -\$\$ & 0  & -\frac{b(\xvec)}{\eps v_\perp}  \\
0 & \frac{\cos(\theta)}{v_\perp} & -\frac{\sin(\theta)}{v_\perp} &-\$ & \frac{b(\xvec)}{\eps v_\perp} & 0\\
\end{array}\right),
\end{gather}
where
\begin{gather} 
\$ = \frac{v_\parallel}{\eps} (\fracp{b}{x_2}(\xvec) +\fracp{b}{x_3}(\xvec)) \text{ ~ and ~ }
\$\$ = \frac{v_\parallel}{\eps v_\perp} (\sin(\theta)\fracp{b}{x_2}(\xvec) +\frac{1+\sin^2(\theta)}{\cos(\theta)}\fracp{b}{x_3}(\xvec)).
\end{gather}

In more general cases, we do not give the Poisson Matrix expression. Yet we mention a very important fact which is that
\begin{gather}
\label{CyCooIVel795} 
\ds (\widetilde \Pcal^\eps_{\!\eta}(\xvec,v_\parallel,v_\perp,\theta))_{56} = 
(\widetilde \Pcal^\eps_{\!\eta}(\xvec,v_\parallel,v_\perp,\theta))_{v_\perp\theta} = \frac{|\Bvec(\xvec)|}{\eps v_\perp}> 0, 
\end{gather}
is always true.
This is Important for the Darboux Method which will be reached in a few lines. We will denote:
\begin{gather}
\label{CyCooIVel797} 
\frac{|\Bvec(\xvec)|}{\eps v_\perp} = \omega(\xvec,v_\perp).
\end{gather}
%
%
\section{Darboux Algorithm}

\paragraph{Darboux Algorithm Target -} The third step is the application of a mathematical algorithm, so called the Darboux Algorithm,
to build a Coordinate System $(\yvec,u_\parallel,k,\theta)$ in which the Poisson Matrix has the required form, given by  \eqref{HamSystSh1},
to apply the Key Result. In fact, in order to manage the small parameter $\eps$, we will build the
Coordinate System $(\yvec,u_\parallel,k,\theta)$ in order to get  $\overline\Pcal^\eps_{\!\eta}$  with the following form:
\begin{gather}
\label{Dab340} 
                                  \left(\begin{array}{c|cc}
                    \text{\huge${\cal M}$} & \begin{array}{c}0 \vspace{-3pt}\\ : \\0 \end{array} & \begin{array}{c}0 \vspace{-3pt}\\ : \\0 \end{array} \\
                    \hline
                     0 \cdots 0& 0 & \frac1\eps\\
                     0 \cdots 0& - \frac1\eps  & 0\\
                                   \end{array}\right).
\end{gather}
We introduce the following notations to manage change-of-coordinates mappings:
$(\yvec,u_\parallel,$ $k, \theta) = \ytrf(\xvec,v_\parallel,v_\perp,\theta)$ and $
(\xvec,v_\parallel,v_\perp,\theta) = \xtrf(\yvec,u_\parallel,k,\theta)$, ${(\xtrf=\ytrf^{-1})}$.\\
An important and constitutive fact in the Darboux Algorithm is that the $\theta$-variable is left unchanged.\\

Now, Since
\begin{gather}
\label{Dab345} 
(\overline\Pcal^\eps_{\!\eta}(\yvec,v_\parallel,k,\theta))_{ij} = \poibrack{\ytrf_{\!i}}{\ytrf_{\!j}}(\xtrf(\yvec,v_\parallel,k,\theta)),
\poibrack{\ytrf_{\!i}}{\ytrf_{\!j}}= (\nabla{\ytrf_{\!i}}) \cdot (\widetilde\Pcal^\eps (\nabla{\ytrf_{\!j}})),
\end{gather}
the bottom-right of form given in \eqref{Dab340}, results from: 
\begin{gather}
\label{Dab348} 
\poibrack{\ytrf_{\!6}}{\ytrf_{\!5}}=-\frac1\eps.
\end{gather}
\begin{remark}
\label{rq:Dab348} 
If we write  $\ytrf_{\!6}=\ytrf_{\!\theta}$ and  $\ytrf_{\!5}=\ytrf_{\!k}$,  \eqref{Dab348} may be also read $\poibrack{\ytrf_{\!\theta}}{\ytrf_{\!k}}=-\frac1\eps$.
In articles of physicists, this last equation reads $\poibrack{\theta}{k}=-\frac1\eps$.
\end{remark}
In the same way, the fact that the two last lines (or columns) contain only zeros results from: 
\begin{gather}
\label{Dab347} 
\begin{aligned}
&\ds \poibrack{\ytrf_{\!1}}{\ytrf_{\!5}}= 0, &
&\ds \poibrack{\ytrf_{\!1}}{\ytrf_{\!6}}= 0,\\
&\ds \poibrack{\ytrf_{\!2}}{\ytrf_{\!5}}= 0, &
&\ds \poibrack{\ytrf_{\!2}}{\ytrf_{\!6}}= 0,\\
&\ds \poibrack{\ytrf_{\!3}}{\ytrf_{\!5}}= 0, &
&\ds \poibrack{\ytrf_{\!3}}{\ytrf_{\!6}}= 0,\\
&\ds \poibrack{\ytrf_{\!4}}{\ytrf_{\!5}}= 0, &
&\ds \poibrack{\ytrf_{\!4}}{\ytrf_{\!6}}= 0.
\end{aligned}
\end{gather}
\begin{remark}
\label{rq:Dab349} 
Using the same conventions as in Remark  \ref{rq:Dab348},  \eqref{Dab347} may also read
\begin{gather}
\label{Dab347-1} 
\begin{aligned}
&\ds \poibrack{\ytrf_{\!y_1}}{\ytrf_{\!k}} =0, \text{ or } \poibrack{y_1}{k} = 0, ~~~~&
&\ds \poibrack{\ytrf_{\!y_1}}{\!\ytrf_{\!\theta}} =0, \text{ or } \poibrack{y_1}{\theta}= 0,\\
&\ds \poibrack{\ytrf_{\!y_2}}{\!\ytrf_{\!k}}=0, \text{ or } \poibrack{y_2}{k}= 0, &
&\ds \poibrack{\ytrf_{\!y_2}}{\!\ytrf_{\!\theta}}=0, \text{ or } \poibrack{y_2}{\theta}= 0,\\
&\ds \poibrack{\ytrf_{\!y_3}}{\!\ytrf_{\!k}}=0, \text{ or } \poibrack{y_3}{k}= 0, &
&\ds \poibrack{\ytrf_{\!y_3}}{\!\ytrf_{\!\theta}}=0, \text{ or } \poibrack{y_3}{\theta}= 0,\\
&\ds\poibrack{\ytrf_{\!y_4}}{\!\ytrf_{\!k}}=0, \text{ or } \poibrack{v_\parallel}{k}\!= 0, &
&\ds \poibrack{\ytrf_{\!v_\parallel}}{\!\ytrf_{\!\theta}}=0, \text{ or } \poibrack{v_\parallel}{\theta}\!= 0.
\end{aligned}
\end{gather}
\end{remark}

Equations  \eqref{Dab348} and \eqref{Dab347} are hyperbolic PDEs that need to be solve to get change-of-coordinates mapping $\ytrf$.
\paragraph{First equation processing -}
Since (and this is a consequence of the fact that the $\theta$-variable is left unchanged in the sought change-of-coordinates)
\begin{gather}
\label{Dab349} 
\nabla\ytrf_{\!6} (=\nabla\ytrf_{\theta}) = (0,0,0,0,0,1)^T, 
\end{gather}
we deduce $\poibrack{\!\ytrf_{\!6}}{\ytrf_{\!5}} = (\nabla{\ytrf_{\!6}}) \cdot (\widetilde\Pcal^\eps (\nabla{\ytrf_{\!5}}))$ is the last component of 
$(\widetilde\Pcal^\eps (\nabla{\ytrf_{\!5}}))$. Hence, equation \eqref{Dab347} reads 
\begin{gather}
\label{Dab350} 
 F_1\fracp{\ytrf_{\!5}}{x_1} + F_2 \fracp{\ytrf_{\!5}}{x_2} +F_3 \fracp{\ytrf_{\!5}}{x_3} +F_\parallel  \fracp{\ytrf_{\!5}}{v_\parallel} 
+ \omega \fracp{\ytrf_{\!5}}{v_\perp}  = -\frac1\eps,
\end{gather}
where $F_n$ and $\omega$ are functions of $(\xvec,v_\parallel,v_\perp,\theta)$. Function $\omega$ is given by \eqref{CyCooIVel797}
and is positive. Solving \eqref{Dab350} will give $\ytrf_{\!5}(\xvec,v_\parallel,v_\perp,\theta)$ which is the expression of component $k$
of the Darboux Coordinates in terms of the Cylindrical in Velocity Coordinates $(\xvec,v_\parallel,v_\perp,\theta)$. 

To solve this equation we will use the Method of Characteristics.
%

\paragraph{Method of Characteristics -}
Dividing by $\omega$, which is possible since $\omega$ is positive, \eqref{Dab350} gives
\begin{align}
\label{Dab356} 
& \fracp{\ytrf_{\!5}}{v_\perp} + \eps\frac{v_\perp F_1}{|\Bvec|}\fracp{\ytrf_{\!5}}{x_1} + \eps\frac{v_\perp F_2}{|\Bvec|}\fracp{\ytrf_{\!5}}{x_2} 
+\eps\frac{v_\perp F_3}{|\Bvec|}\fracp{\ytrf_{\!5}}{x_3}
+\eps\frac{v_\perp F_\parallel }{|\Bvec|} \fracp{\ytrf_{\!5}}{v_\parallel}  = \frac{v_\perp}{|\Bvec|}.
\end{align}
In order to get a solution to this equation, we need to add a boundary (or initial) condition in a point where the involved vector field 
$(\eps v_\perp/|\Bvec|)(F_1, F_2, F_3,F_\parallel)$ does not vanish.  Hence we set
\begin{align}
\label{Dab357} 
& {\ytrf_{\!5}}_{|v_\perp=\nu} =0 \text{ ~~~~~ for a small } \nu>0.
\end{align}

~

We consider the following characteristics $(\Xcalvec,\Vcal_\parallel)= (\Xcal_1,\Xcal_2,\Xcal_3,\Vcal_\parallel)$:
\begin{gather}
\label{Dab359} 
\begin{aligned}
&\Xcal_1(v_\perp; \xvec,v_\parallel, u_\perp )\text{ such that } &
&\ds \fracp{\Xcal_1}{v_\perp} = \eps\frac{v_\perp F_1(\Xcal_1,\Xcal_2,\Xcal_3,\Vcal_\parallel,v_\perp,\theta)}{|\Bvec|(\Xcal_1,\Xcal_2,\Xcal_3)},&
&\Xcal_1(u_\perp)\!=\!x_1,
\\
&\Xcal_2(v_\perp; \xvec,v_\parallel, u_\perp )\text{ such that } &
&\ds \fracp{\Xcal_2}{v_\perp} = 
\eps\frac{v_\perp F_2(\Xcal_1,\Xcal_2,\Xcal_3,\Vcal_\parallel,v_\perp,\theta)}{|\Bvec|(\Xcal_1,\Xcal_2,\Xcal_3)},&
&\Xcal_2(u_\perp)\!=\!x_2,
\\
&\Xcal_3(v_\perp; \xvec,v_\parallel, u_\perp )\text{ such that } &
&\ds \fracp{\Xcal_3}{v_\perp} = \eps\frac{v_\perp F_3(\Xcal_1,\Xcal_2,\Xcal_3,\Vcal_\parallel,v_\perp,\theta)}{|\Bvec|(\Xcal_1,\Xcal_2,\Xcal_3)},&
&\Xcal_3(u_\perp)\!=\!x_3,
\\
&\Vcal_\parallel(v_\perp; \xvec,v_\parallel, u_\perp )\text{ such that } &
&\ds \fracp{\Vcal_\parallel}{v_\perp} 
= \eps\frac{v_\perp F_\parallel(\Xcal_1,\Xcal_2,\Xcal_3,\Vcal_\parallel,v_\perp,\theta)}{|\Bvec(\Xcal_1,\Xcal_2,\Xcal_3)},&
&\Vcal_\parallel(u_\perp)\!=\!v_\parallel,\\
\end{aligned}
\end{gather}
and the solution to \eqref{Dab356} - \eqref{Dab357} is given by
\begin{multline}
\label{Dab363} 
 \ytrf_{\!5} (\xvec,v_\parallel,v_\perp,\theta) = 
 \ytrf_{\!5} (\Xcalvec(\nu; \xvec,v_\parallel, v_\perp),\Vcal_\parallel(\nu; \xvec,v_\parallel, v_\perp),\nu,\theta) +
 \int_\nu^{v_\perp} \frac{s}{|\Bvec(\Xcalvec(s; \xvec,v_\parallel, v_\perp))|} ds\\
 = \int_\nu^{v_\perp} \frac{s}{|\Bvec(\Xcalvec(s; \xvec,v_\parallel, v_\perp))|} ds.
 \end{multline}
(The last equality is gotten because of \eqref{Dab357}.)

Rewriting system \eqref{Dab359} in the following compact form,
\begin{gather}
\label{Dab366} 
\fracp{\ds \begin{pmatrix}\Xcalvec \\ \Vcal_\parallel \end{pmatrix}}{v_\perp} = \eps \Fvec (\Xcalvec_1,\Vcal_\parallel,v_\perp,\theta),
\begin{pmatrix}\Xcalvec \\ \Vcal_\parallel \end{pmatrix}(u_\perp) =\begin{pmatrix} \xvec \\ v_\parallel \end{pmatrix},
\end{gather}
with $\Fvec = (v_\perp/|\Bvec|)(F_1, F_2, F_3,F_\parallel)$ and applying the formula that gives the expansion of the solution
of a dynamical system in terms of its parameter, we get that 

\begin{gather}
\label{Dab3366} 
\ds \begin{pmatrix}\Xcalvec \\ \Vcal_\parallel \end{pmatrix}=
\ds \begin{pmatrix}\Xcalvec \\ \Vcal_\parallel \end{pmatrix}(v_\perp; \xvec,v_\parallel, u_\perp),
\end{gather}
writes 
\begin{multline}
\label{Dab367} 
\begin{pmatrix}\Xcalvec \\ \Vcal_\parallel \end{pmatrix} 
= \begin{pmatrix} \xvec \\ v_\parallel \end{pmatrix}
+ \frac{v_\perp}{1!} \eps\Fvec (\xvec,v_\parallel,v_\perp,\theta) + \frac{v_\perp^2}{2!}  \Lcal_{\eps\Fvec}(\eps\Fvec) (\xvec_1,v_\parallel,v_\perp,\theta) 
\\
~ \hspace{5cm}+ \frac{v_\perp^3}{3!}  \Lcal_{\eps\Fvec}^2(\eps\Fvec) (\xvec_1,v_\parallel,v_\perp,\theta) + \dots
\\
= \begin{pmatrix} \xvec \\ v_\parallel \end{pmatrix}
+ \frac{\eps v_\perp}{1!} \Fvec (\xvec,v_\parallel,v_\perp,\theta) + \frac{\eps^2 v_\perp^2}{2!}  \Lcal_{\Fvec}(\Fvec) (\xvec_1,v_\parallel,v_\perp,\theta) 
\\
+ \frac{\eps^3 v_\perp^3}{3!}  \Lcal_{\Fvec}^2(\Fvec) (\xvec_1,v_\parallel,v_\perp,\theta) + \dots
\end{multline}
where $\Lcal_\Fvec(\Fvec)$  is the the Lie Derivative of vector field $\Fvec$ in the  direction of vector field $\Fvec$, 
$\Lcal^2_\Fvec(\Fvec)= \Lcal_\Fvec(\Lcal_\Fvec(\Fvec))$  and $\Lcal^3_\Fvec(\Fvec)= \Lcal_\Fvec(\Lcal^2_\Fvec(\Fvec))$, and so on.

Rewriting this expansion, we get
\begin{gather}
\label{Dab365} 
\Xcalvec = \xvec + \eps v_\perp \Xcalvec^1 + \eps^2v_\perp^2  \Xcalvec^2 + \eps^3v_\perp^3  \Xcalvec^3 + \dots
\\
\label{Dab366.1} 
\Vcal_\parallel =v_\parallel + \eps v_\perp \Vcal_\parallel^1 + \eps^2 v_\perp^2  \Vcal_\parallel^2 + \eps^3 v_\perp^3  \Vcal_\parallel^3+ \dots
\end{gather}
with
\begin{gather}
\label{Dab368} 
\begin{pmatrix}\Xcalvec^1 \\ \Vcal_\parallel^1\end{pmatrix} = \Fvec (\xvec,v_\parallel,v_\perp,\theta),
\begin{pmatrix}\Xcalvec^2 \\ \Vcal_\parallel^2\end{pmatrix} =\frac{1}{2!} \Lcal_{\Fvec}(\Fvec) (\xvec,v_\parallel,v_\perp,\theta),
\begin{pmatrix}\Xcalvec^3 \\ \Vcal_\parallel^3\end{pmatrix} =\frac{1}{3!} \Lcal^2_{\Fvec}(\Fvec) (\xvec,v_\parallel,v_\perp,\theta),
\dots\end{gather}

Now, injecting \eqref{Dab366.1} in expression \eqref{Dab363} on $\ytrf_{\!5}$, we get:
 \begin{multline}
\label{Dab3368} 
 \ytrf_{\!5} (\xvec,v_\parallel,v_\perp,\theta) =
  \int_\nu^{v_\perp} \frac{s}{|\Bvec(\Xcalvec(s; \xvec,v_\parallel, v_\perp))|} ds =\\
  \int_\nu^{v_\perp} \frac{s}{|\Bvec(\xvec)|} ds + \eps \int_\nu^{v_\perp} s^2 \; \Tcal^1(\frac{1}{|\Bvec(\xvec)|})\Xcalvec^1 ds + 
  \\ ~~~~~~~~~~~~~~~~~~~~~~~~~~~~~~~~~~~
  \eps^2 \int_\nu^{v_\perp} s^3 \; \Big(\Tcal^2(\frac{1}{|\Bvec(\xvec)|})\Xcalvec^1 + \Tcal^1(\frac{1}{|\Bvec(\xvec)|})\Xcalvec^2\Big) ds + \dots
  \\
 = \frac{(v_\perp-\nu)^2}{2|\Bvec(\xvec)|} + \eps \int_\nu^{v_\perp} s^2 \; \Tcal^1(\frac{1}{|\Bvec(\xvec)|})\Xcalvec^1 ds + 
~~~~~~~~~~~~~~~~ \\ 
  \eps^2 \int_\nu^{v_\perp} s^3 \; \Big(\Tcal^2(\frac{1}{|\Bvec(\xvec)|})\Xcalvec^1 + \Tcal^1(\frac{1}{|\Bvec(\xvec)|})\Xcalvec^2\Big) ds +  \dots,
 \end{multline}
where $\Tcal^i(1/|\Bvec(\xvec)|)$ is the $i^\text{th} $ coefficient of the  Taylor expansion of $1/|\Bvec(\xvec)|$.

~

Formula \eqref{Dab3368} gives the expression of  $\ytrf_{\!5}$, i.e. the expression of new variable $k$ in terms of 
Cynlindrical Coordinates $(\xvec,v_\parallel,v_\perp,\theta)$, as an expansion in $\eps$.
\begin{remark}
In \eqref{Dab3368}, it is possible to choose $\nu$ as small as we want. Hence, the first term of the expansion is essentially
the Magnetic Moment $v_\perp/|\Bvec(\xvec)|$.
\end{remark}

\paragraph{On other equations and on Poisson Matrix in Darboux Coordinates}
Equation \eqref{Dab348} was processed and gave expression of $k$.
Equations \eqref{Dab347} can be processed using similar Methods of Characteristics.  A special attention needs to be given
to the fact that $\ytrf_{\!1}$, \dots $\ytrf_{\!4}$ are solutions of a PDE involving $\ytrf_{\!5}$ and of another one involving
$\ytrf_{\!6}$ that need to be tackled together.  
They will give $\yvec$ and $u_\parallel$ as expansions in $\eps$. Once this done, we will have
\begin{gather}
\label{upsexp} 
\ytrf = \ytrf_0 + \eps\ytrf_1 + \eps^2\ytrf_2 + \eps^3\ytrf_3 + \dots,
\end{gather}
or, in other words, the expression of $(\yvec,u_\parallel,k)$ in terms of $(\xvec,v_\parallel,v_\perp,\theta)$ and as an expansion in $\eps$.

~
 
The other terms of the  new Poisson matrix $\overline \Pcal^\eps_{\!\eta}(\yvec,u_\parallel,k,\theta)$ are given computing:
$\ds(\widehat \Pcal^\eps_{\!\eta})_{12} =  \poibrack{\ytrf_{\!1}}{\!\ytrf_{\!2}}(=\!\!\poibrack{\ytrf_{\!y_1}}{\!\ytrf_{\!y_2}}\!\!=\!\!\poibrack{y_1}{y_2})$,
$\ds(\widehat \Pcal^\eps_{\!\eta})_{13} =  \poibrack{\ytrf_{\!1}}{\!\ytrf_{\!3}}(=\!\!\poibrack{\ytrf_{\!y_1}}{\!\ytrf_{\!y_3}}\!\!=\!\!\poibrack{y_1}{y_3})$,
\dots.

\paragraph{Hamiltonian Function in Darboux Coordinates -}
We now have to compute the Hamiltonian Function in the Darboux Coordinates. For this, we build from expansion \eqref{upsexp}
of $\ytrf$, an asymptotic expansion of $\xtrf = \ytrf^{-1}$:
\begin{gather}
\label{xiexp} 
\xtrf = \xtrf_0 + \eps\xtrf_1 + \eps^2\xtrf_2 + \eps^3\xtrf_3 + \dots.
\end{gather}
Then, we use the expression of the Hamiltonian Function in Cylindrical in Velocity Coordinates:
\begin{gather}
\ds \widetilde H^\eps_{\!\eta}\!(\xvec,v_\parallel,v_\perp,\theta)= 
\frac 12 \big(v_\parallel ^2+v_\perp^2\big) + \Phi_0(\xvec) + {\eta} \Phi_1(\frac{\xvec}{\eta}),
\end{gather}
and the following formula, which is gotten from  \eqref{PmPbCcf28},
\begin{multline}
\label{Dab5368} 
\ds \overline H^\eps_{\!\eta}\!(\yvec,u_\parallel,k,\theta)= \ds \widetilde H^\eps_{\!\eta}\!(\xtrf(\yvec,u_\parallel,k,\theta))
= \widetilde H^\eps_{\!\eta}\!( \xtrf_0 + \eps\xtrf_1 + \eps^2\xtrf_2 + \eps^3\xtrf_3 + \dots) =
\\
 \widetilde H^\eps_{\!\eta}\!( \xtrf_0) + \eps \Tcal^1(\widetilde H^\eps_{\!\eta}\!)( \xtrf_0)\cdot \xtrf_1+ \dots=
 \\
 \ \frac{u_\parallel^2}2 +|\Bvec(\yvec)|k + \Phi_0(\yvec) 
+ {\eta} \Phi_1(\frac{\yvec}{\eta}) + 
 \eps  \overline H_{\!1,\eta}(\yvec,u_\parallel,k,\theta) + \eps^2  \overline H_{\!2,\eta}(\yvec,u_\parallel,k,\theta)+ \dots.
 \end{multline}
Notice that for this expansion to be valid, it is necessary to have the property we mention on page \pageref{gghh}: 
$\eta= \eps^{1-\kappa}$ for $\kappa>0$. \\

In expression  \eqref{Dab5368}, there is an important fact for the setting out of the to come Lie  
Transform based Method: the first term is independent of $\theta$.
%
%
%
\section{Lie Transform based Method}
%
%

\paragraph{Lie Transform based Method Target -}
As a result of the Darboux Algorithm, we obtained a Poisson Matrix 
$\overline \Pcal^\eps_{\!\eta}(\yvec,u_\parallel,k,\theta)$ with the required form to apply the Key Result (Theorem  \ref{KrTh}), but
 the resulting Hamiltonian Function
\begin{gather}
\label{Lie3401} 
\ds \overline H^\eps_{\!\eta}\!(\yvec,u_\parallel,k,\theta)= \overline H_{\!0,\eta}(\yvec,u_\parallel,k) 
+ \eps \overline H_{\!1,\eta}(\yvec,u_\parallel,k,\theta)+  \eps^2 \overline H_{\!2,\eta}(\yvec,u_\parallel,k,\theta)+\dots,
\end{gather}
depends on $\theta$. And we need to make this dependency to vanish.

Let us notice two important  facts, which are in fact  linked:
The first term in the asymptotic expansion  \eqref{Dab5368} does not depend on $\theta$ and the Hamiltonian
Function expressed in the Cylindrical in Velocity Coordinates does not depend on $\theta$. The consequence
of those facts is that it is certainly possible to make the last variable to vanish using a mapping parametrized by  $\eps$  
and close to identity for small  $\eps$.
\\
Moreover, we need to build this sought mapping in such a way that it does not change the 
Poisson Matrix expression. This means that, as viewed as functions, $\overline\Pcal^\eps_{\!\eta}$ and $\widehat \Pcal^\eps_{\!\eta}$
must be the same, i.e.:
\begin{gather}
\ds \widehat \Pcal^\eps_{\!\eta}(\zvec,w_\parallel,j,\gamma)=\overline \Pcal^\eps_{\!\eta}(\zvec,w_\parallel,j,\gamma)
\text{ or }
\overline \Pcal^\eps_{\!\eta}(\yvec,u_\parallel,k,\theta) = \widehat \Pcal^\eps_{\!\eta}(\yvec,u_\parallel,k,\theta),
\end{gather}
for any $\zvec,w_\parallel,j,\gamma$ or $\yvec,u_\parallel,k,\theta$.
Yet, among changes of variables, the symplectic ones do not change the Poisson Matrix expression.
And among  symplectic changes of variables,  are the flows of Hamiltonian Vector Fields, which are moreover  close to identity
for small values of their parameter.\\

Hence, the Lie Transform based Method consists in building a change of variables, parametrized by $\eps$, reading
\begin{gather}
\label{Lie3402} 
(\yvec,u_\parallel,k,\theta) \mapsto (\zvec,w_\parallel,j,\gamma) = \zetavec(\eps;\yvec,u_\parallel,k,\theta),
\end{gather}
with
\begin{gather}
\label{Lie3404} 
\zetavec(\eps;\yvec,u_\parallel,k,\theta) = (\yvec,u_\parallel,k,\theta) + \eps \zetavec_1(\yvec,u_\parallel,k,\theta),
+\eps^2 \zetavec_2(\yvec,u_\parallel,k,\theta)+ \dots,
\end{gather}
 which is the flow of parameter $\eps$ of a Hamiltonian Vector Fields.
\paragraph{The way to do -}
Looking for $\zetavec$ as the flow of a Hamiltonian Vector Field means looking for a  Hamiltonian Function:
\begin{gather}
\label{Lie3405} 
G(\eps;\yvec,u_\parallel,k,\theta) = G_0(\yvec,u_\parallel,k,\theta) + \eps G_1(\yvec,u_\parallel,k,\theta)
+ \eps^2 G_2(\yvec,u_\parallel,k,\theta) + \dots,
\end{gather}
such that $ \zetavec(\eps;\yvec,u_\parallel,k,\theta)$ is  solution to:
\begin{gather}
\label{Lie3406} 
\fracp\zetavec\eps = \overline\Pcal^\eps_{\!\eta} \nabla G, ~~ \zetavec(\eps=0;\yvec,u_\parallel,k,\theta ) = (\yvec,u_\parallel,k,\theta).
\end{gather}

Hence the target becomes to find functions  $G_0$, $G_1$,  $G_2$, \dots,  such that:
\begin{gather}
\label{Lie3408} 
\widehat H^\eps_{\!\eta}\!(\zvec,w_\parallel,j) =\overline H^\eps_{\!\eta}\!(\ytrfb(\zvec,w_\parallel,j,\gamma)) \text{ where }\ytrfb=\ztrf^{-1}.
\end{gather}

%

\paragraph{A result -}
We will use the following theorem
\begin{theorem}
For any two Hamiltonian Functions:
\begin{gather}
\label{Lie3410} 
\ds \overline H^\eps_{\!\eta}\!(\yvec,u_\parallel,k,\theta)= \overline H_{\!0,\eta}(\yvec,u_\parallel,k) 
+ \eps \overline H_{\!1,\eta}(\yvec,u_\parallel,k,\theta)+  \eps^2 \overline H_{\!2,\eta}(\yvec,u_\parallel,k,\theta)+\dots,
\end{gather}
and
\begin{gather}
\label{Lie3411} 
\ds \widehat H^\eps_{\!\eta}\!(\zvec,w_\parallel,j)= \widehat H_{\!0,\eta}(\zvec,w_\parallel,j) 
+ \eps \widehat H_{\!1,\eta}(\zvec,w_\parallel,j)+  \eps^2 \widehat H_{\!2,\eta}(\zvec,w_\parallel,j)+\dots,
\end{gather}
expressed in two different variable systems ($\yvec,u_\parallel,k,\theta$ and $\zvec,w_\parallel,j,\gamma$) , with the property that 
\begin{gather}
\widehat H_{\!0,\eta}(\zvec,w_\parallel,j) = \overline H_{\!0,\eta}(\zvec,w_\parallel,j),
\end{gather}
There exists a Hamitonian Function writing
\begin{gather}
\label{Lie3413} 
G(\eps;\yvec,u_\parallel,k,\theta) = G_0(\yvec,u_\parallel,k,\theta) + \eps G_1(\yvec,u_\parallel,k,\theta) 
+ \eps^2 G_2(\yvec,u_\parallel,k,\theta) + \dots,
\end{gather}
such that the solution $\zetavec$ of the Hamiltonian Dynamical System associated with $G$ and with parameter $\eps$, i.e. solution to:
\begin{gather}
\label{Lie3414} 
\fracp\zetavec\eps = \overline\Pcal^\eps_{\!\eta} \nabla G, ~~ \zetavec(\eps=0;\zvec,w_\parallel,j,\gamma) = (\yvec,u_\parallel,k,\theta),
\end{gather}
is such that 
\begin{gather}
\ds \widehat H^\eps_{\!\eta}\!(\zvec,w_\parallel,j) = \overline H^\eps_{\!\eta}\!(\ytrfb(\zvec,w_\parallel,j,\gamma))  \text{ where } \ytrfb=\ztrf^{-1}.
\end{gather}
Moreover, the Hamiltonian Function $g(\eps;\zvec,w_\parallel,j,\gamma)$ such that $\ytrfb$ is solution to
\begin{gather}
\label{Lie7414} 
\fracp\ytrfb\eps = \overline\Pcal^\eps_{\!\eta} \nabla g, ~~ \ytrfb(\eps=0;\yvec,u_\parallel,k,\theta) = (\zvec,w_\parallel,j,\gamma),
\end{gather}
writes
\begin{gather}
\label{Lie3416} 
g(\eps;\zvec,w_\parallel,j,\gamma) = g_0(\zvec,w_\parallel,j,\gamma) +\eps g_1(\zvec,w_\parallel,j,\gamma) 
+ \eps^2 g_2(\zvec,w_\parallel,j,\gamma) + \dots,
\end{gather}
where $g_0$, $g_1$,  $g_2$, \dots are given by 
\begin{multline}
\label{Lie3418} 
\poibrack{g_0}{\overline H_{\!0,\eta}} = \Ocal_0(\overline H_{\!0,\eta}), ~~
\poibrack{g_1}{\overline H_{\!0,\eta}} = \Ocal_1(\overline H_{\!0,\eta},\overline H_{\!1,\eta},\widehat H_{\!1,\eta},g_0), ~~\\ 
\poibrack{g_2}{\overline H_{\!0,\eta}} = 
\Ocal_2(\overline H_{\!0,\eta},\overline H_{\!1,\eta},\widehat H_{\!1,\eta},\overline H_{\!2,\eta},\widehat H_{\!2,\eta}, g_0,g_1), \dots.
\end{multline}
for differential operators $\Ocal_0$, $\Ocal_1$, $\Ocal_2$, \dots defined by recursive formula.
\end{theorem}

\begin{remark}
Despite we do not explicitly write the dependency of $\ytrfb$,  $g_0$, $g_1$,  $g_2$, \dots, $\ytrfb$ and $G_0$, $G_1$,  $G_2$, \dots
with respect $\eta$ to all those functions do depend on $\eta$. Here again, to be able to write the expansions with respect to $\eps$, it is 
necessary to have the property mentioned on page \pageref{gghh}: 
$\eta= \eps^{1-\kappa}$ for $\kappa>0$. \\
\end{remark}

%
%

\paragraph{The Lie Transform based Method - } As in the case of the Darboux Algorithm, we do not do the computations that give 
$\ztrf$ and $\ytrfb$. We only give the steps that  allows us to get them. The Lie Transform based Method may be summarized as: 
\begin{enumerate}
\item Fix : $\widehat H_{\!1,\eta}(\zvec,w_\parallel,j)$, $\widehat H_{\!2,\eta}(\zvec,w_\parallel,j), \dots$. \\
They can be fixed, a priory, with no restriction. Nevertheless those functions are involved in the PDEs of  \eqref{Lie3418}. Hence,  they
need to be chosen in a way that leads to PDEs which are as simple as possible to solve.
\item Solve PDE in  \eqref{Lie3418} to get $g_0$, $g_1$,  $g_2$ recursively.\\
They are quite complex but we can solve them using the Method of  Characteristics.
\item Get  $\ztrf$ solving  \eqref{Lie3414}.\\
It is gotten as an expansion in $\eps$.
\item Compute expression of the Lie Variables: $ (\zvec,w_\parallel,j,\gamma) = \zetavec(\yvec,u_\parallel,k,\theta)$
\item Since $\widehat H_{\!1,\eta}$, $\widehat H_{\!2,\eta}$, \dots were fixed in the begining of the process, 
it is not necessary to compute $\ytrfb$ to get: 
$\ds \widehat H^\eps_{\!\eta}\!(\zvec,w_\parallel,j) =\overline H^\eps_{\!\eta}\!(\ytrfb(\zvec,w_\parallel,j,\gamma))$. Nonetheless, $\ytrfb$ may
be gotten as an expansion in $\eps$ by inverting expansion of $\ztrf$.
\end{enumerate}
\section{Brief conclusion}
At the end of the day, in the Lie Coordinates System, the trajectory $(\Zvec,W_\parallel,J,\Gamma)$ of a charged particle is solution to
a Dynamical System which has the following form:
\begin{align}
\label{Lie5810} 
&\fracp{\Gamma}{t} = \text{Something complicated},
\\
\label{Lie5811} 
&\fracp{J}{t} =  0,
\\
\label{Lie5813} 
&\fracp{\ds \begin{pmatrix} \Zvec \\  W_\parallel \end{pmatrix}}{t} = \text{Something independent of } \Gamma.
\end{align}
This form is essentially the same as the one of \eqref{GyrMod}.

~

As a conclusion, we just say that the method which is summarized here is the one that allows us to get  \eqref{GyrMod} or 
any other Gyro-Kinetic Approximation Model.

~

Insisting one more time, the important facts in system  \eqref{Lie5810} -  \eqref{Lie5813} is that the evolution of the components 
$\Zvec$ and $W_\parallel$ of the trajectory
and of its $\gamma$-component $\Gamma$ are uncoupled  and that component $J$ does not evolve and then becomes a parameter.
Hence, even if the $\gamma$-component $\Gamma$ is not computed, it does not preclude the computation of  $J$, $\Zvec$ and 
$W_\parallel$ solving a collection, parametrized by  $J$,  of four-dimensional Dynamical Systems.
%
%
%
%
%
\bibliographystyle{plain}
\bibliography{biblio}

\end{document}